\documentclass[11pt,reqno]{amsart}
\usepackage{amsfonts}
\usepackage{mathrsfs}
\usepackage{graphicx} 
\usepackage{epsfig}
\usepackage{amsmath, amsthm}
\usepackage{amssymb}

\numberwithin{equation}{section}
\date{}

\theoremstyle{plain}

\theoremstyle{definition}

\numberwithin{equation}{section}

\theoremstyle{plain}
\newtheorem{thm}{Theorem}[section]
\newtheorem{prop}[thm]{Proposition}
\newtheorem{lem}[thm]{Lemma}
\newtheorem{cor}[thm]{Corollary}

\theoremstyle{definition}

\newtheorem{rem}[thm]{Remark}
\newtheorem{defn}[thm]{Definition}
\newtheorem{eg}[thm]{Example}
\newtheorem{subtitle}[thm]{}
\newtheorem{ex}{Exercise}[section]
\numberwithin{equation}{section}

\def\a{\alpha}
\def\b{\beta}
\def\d{\delta}
\def\D{\triangle}
\def\e{\epsilon}
\def\g{\gamma}
\def\G{\Gamma}
\def\K{\nabla}
\def\l{\lambda}

\def\n{\,\vert\,}

\def\cb{{\mathcal{B}}}

\def\cd{{\mathcal{D}}}

\def\cg{{\mathcal{G}}}

\def\ci{{\mathcal{I}}}

\def\cl{{\mathcal{L}}}
\def\cm{{\mathcal{M}}}
\def\cn{{\mathcal{N}}}

\def\ct{{\mathcal{T}}}

\def\li{\langle}
\def\ri{\rangle}
\def\n{\ \vert\ }
\def\tr{{\rm tr}}
\def\bs{\bigskip}
\def\ms{\medskip}

\def\ni{\noindent}
\def\ti{\tilde}
\def\p{\partial}

\def\Im{{\rm Im\/}}
\def\I{{\rm I\/}}

\def\diag{{\rm diag}}

\def\A{\mathbb{A}}
\def\C{\mathbb{C}}

\def\R{\mathbb{R} }
\def\Z{\mathbb{Z}}

\newcommand{\beq}{\begin{equation}}
\newcommand{\eeq}{\end{equation}}
\newcommand{\beg}{\begin{eg}}
\newcommand{\eeg}{\end{eg}}
\newcommand{\bthm}{\begin{thm}}
\newcommand{\ethm}{\end{thm}}
\newcommand{\bprop}{\begin{prop}}
\newcommand{\eprop}{\end{prop}}
\newcommand{\bcor}{\begin{cor}}
\newcommand{\ecor}{\end{cor}}
\newcommand{\blem}{\begin{lem}}
\newcommand{\elem}{\end{lem}}
\newcommand{\bca}{\begin{cases}}
\newcommand{\eca}{\end{cases}}
\newcommand{\brem}{\begin{rem}}
\newcommand{\erem}{\end{rem}}
\newcommand{\bpm}{\begin{pmatrix}}
\newcommand{\epm}{\end{pmatrix}}
\newcommand{\bbm}{\begin{bmatrix}}
\newcommand{\ebm}{\end{bmatrix}}
\newcommand{\bvm}{\begin{vmatrix}}
\newcommand{\evm}{\end{vmatrix}}
\newcommand{\bdefn}{\begin{defn}}
\newcommand{\edefn}{\end{defn}}
\newcommand{\bsub}{\begin{subtitle}}
\newcommand{\esub}{\end{subtitle}}
\newcommand{\bex}{\begin{ex}}
\newcommand{\eex}{\end{ex}}
\newcommand{\ben}{\begin{enumerate}}
\newcommand{\een}{\end{enumerate}}

\def\sp{{\rm span \/ }}
\def\rd{{\rm \/ d\/}}

\def\sech{{\rm sech\/}}

\def\anone{\hat A^{(1)}_{2n}}
\def\an2{\hat A_{2n}^{(2)}}
\def\bn1{\hat B_n^{(1)}}
\def\Ker{{\rm Ker\/}}

\def\B{\mathbb{B}}
\def\bbn1{\hat \B_n^{(1)}}
\def\A{\mathbb{A}}
\def\ban2{\hat \A_{2n}^{(2)}}

\def\bh{\backslash}
\def\r0{\R^n \backslash \{0\}}
\def\bc{{\bf c}}
\def\bu{\bullet}

\begin{document}

\title
{Isotropic curve flows on $\R^{n+1, n}$} \today

\author{Chuu-Lian Terng$^\dag$}\thanks{$^\dag$Research supported
in  part by NSF Grant DMS-1109342}
\address{Department of Mathematics\\
University of California at Irvine, Irvine, CA 92697-3875.  Email: cterng@math.uci.edu}
\author{Zhiwei Wu$^*$}\thanks{$^*$Research supported in part by NSF of China under Grant No. 11401327\/}
\address{Department of Mathematics\\ Ningbo University\\ Ningbo, Zhejiang, 315211, China. Email: wuzhiwei@nbu.edu.cn}


\maketitle

\section{Introduction}

Let $\R^{n+1, n}$ be the vector space $\R^{2n+1}$ equipped with the index $n$, non-degenerate bilinear form $\li X, Y \ri=X^tC_nY$, where
\beq\label{am}
C_n=\sum_{i=1}^{2n+1}(-1)^{n+i-1}e_{i, 2n+2-i}.
\eeq
 Let $O(n+1, n)$ denote the group of linear isomorphisms on $\R^{n+1, n}$ preserving $\li \ , \ \ri$. 
A subspace $\ci \subset \R^{n+1, n}$ is called \emph{isotropic} if $\li X, Y \ri =0$ for all  $X, Y \in \ci$.
A maximal isotropic subspace in $\R^{n+1, n}$ has dimension $n$.  

A smooth map $\gamma: \R \rightarrow \R^{n+1,n}$ is called an \emph{isotropic curve} if $\g, \g_s, \ldots$, $ \g_s^{(2n)}$ are linearly independent and the span of
$\gamma, \g_s, \cdots, \g_s^{(n-1)}$
is a maximal isotropic subspace of $\R^{n+1, n}$ for all $s\in \R$.
Note that a curve being isotropic is independent of the choice of parameter. 
It is easy to see that there is an orientation preserving parameter $x$ (unique up to translation) for an isotropic curve such that $\li \g_x^{(n)}, \g_x^{(n)}\ri \equiv 1$. We call such $x$ the {\it isotropic parameter\/} of $\g$. 

Set
$$\cm_{n+1, n}=\{\g:\R\to \R^{n+1,n}\n \g \,\, {\rm isotropic,}\,\, \li \g_x^{(n)}, \g_x^{(n)}\ri \equiv 1\}.$$
We prove that given $\g\in \cm_{n+1, n}$,  there exists a unique smooth map $g:\R\to O(n+1, n)$ such that the $i$-th column is $\g_x^{(i-1)}$ for $1\leq i\leq n+1$ and $g^{-1}g_x$ is of the form 
$$g^{-1}g_x= b+ \sum_{i=1}^n u_i \b_i,$$ where
\beq\label{cj}  b=\sum_{i=1}^{2n}  e_{i+1, i}, \quad \b_i= e_{n+1-i, n+i} + e_{n+2-i, n+1+i}.
\eeq
 We call $g$,  $u_i$, and $u=\sum_{i=1}^n u_i \b_i$ the {\it isotropic moving frame\/}, the {\it $i$-th isotropic curvature\/}, and  the {\it isotropic curvature\/} along $\g$ respectively. 
 
 Set 
\beq\label{aj}
V_n=\oplus_{i=1}^n \R \b_i.
\eeq
 Let $\Psi:\cm_{n+1, n}\to C^\infty(\R, V_n)$ be the {\it isotropy curvature map\/} defined by 
 \beq\label{cd} \Psi(\g)=u =g^{-1}g_x-b= \sum_{i=1}^n u_i \b_i,
 \eeq where $u$ is the isotropy curvature of $\g$. Then $\Psi$ is onto and $\Psi^{-1}(u)$ is an $O(n+1, n)$-orbit for all $u\in C^\infty(\R, V_n)$. Hence $u_1, \ldots, u_n$ form a complete set of {\it differential invariants\/} for  $\g\in \cm_{n+1, n}$.

An {\it isotropic curve flow\/} on $\cm_{n+1, n}$ is an evolution equation on $\cm_{n+1, n}$ of the form
$$\g_t= g \xi(u),$$
where $g(\cdot, t)$, $u(\cdot, t)$ are the isotropic moving frame and isotropic curvature along $\g(\cdot, t)$ respectively and $\xi(u)$ is a $\R^{(2n+1)\times 1}$-valued differential polynomial of $u$ in $x$ variable. Hence isotropic curve flows are invariant under the action of $O(n+1, n)$. 

We note that a map $\g:\R\to \R^{2,1}$ lies in $\cm_{2,1}$ if and only if $\g$ lies in the light cone, is space-like, and is parametrized by its arc-length. The isotropic curvature of such curve is the standard curvature in differential geometry. Hence an isotropic curve flow on $\cm_{2,1}$ is a geometric curve flow for space-like curves on the light cone of $\R^{2,1}$ preserving arc-length. We show that an isotropic curve flow on $\cm_{2,1}$ is of the form 
$$\g_t= (h(q))_x \g - h(q) \g_x$$
for some differential polynomial $h$ of the isotropic curvature $q(\cdot, t)$ of $\g(\cdot, t)$. In particular, 
\begin{align}  &\g_t=q_x\g-q\g_x, \label{bx}\\
& \g_t= -\frac{1}{9} (q_{xxx}-8qq_x)\g +\frac{1}{9} (q_{xx}- 4q^2)\g_x, \label{al}
  \end{align}
are isotropic curve flows. We show that if $\g$ is a solution of \eqref{bx} or \eqref{al}, then the isotropic curvature $q$ is a solution of KdV
 $$q_t= q_{xxx} -3qq_x,$$ or the Kupershmidt-Kaup (KK) equation 
 \beq\label{kk}
q_t=-\frac{1}{9}(q^{(5)}-10qq_{xxx}-25q_xq_{xx}+20q^2q_x)
\eeq
respectively. 
So these two curve flows give natural geometric interpretations of the KdV  and KK equations respectively.

In this paper, we prove the following results:
 \ben
 \item[(a)]  For general $n$, we construct two sequences of isotropic curve flows on $\cm_{n+1, n}$ of B-type and A-type that map to the $\bn1$-KdV and $\an2$-KdV hierarchies under the isotropic curvature map $\Psi$. The third flows of B-type and A-type are
 \begin{align}
\frac{\p \g}{\p t}&=-\frac{3}{n}u_1\g_x+\g_x^{(3)},\label{ar}\\
\frac{\p \g}{\p t}&=-\frac{3}{2n+1}(u_1)_x\g-\frac{6}{2n+1}u_1\g_x+\g_{x}^{(3)}, \label{as}
\end{align}
where
 $u_1(\cdot, t)$ is the first isotropic curvature along  $\g(\cdot, t)$. Note that the flows \eqref{bx} and \eqref{al} are the third and fifth flows of these two sequence when $n=1$.
\item[(b)] Algorithms to compute the Bi-Hamiltonian structures and conservation laws for isotropic curve flows of B- and A-types are given. 
\item[(c)]  We construct B\"acklund transformations and Permutability formulas. 
\een
 
In particular, isotropic curve flows give natural geometric interpretations of the $\bn1$-KdV and $\an2$-KdV flows and  techniques from soliton theory give bi-Hamiltonian, conservation laws, and BTs for these isotropic curve flows. We note that the relation between central curve flows on $\R^n\bh 0$ and the soliton theory of the $\hat A_{n-1}^{(1)}$-KdV hierarchy were considered in \cite{UP95} and \cite{TWa} for $n=2$, in \cite{CIM13} for $n=3$, and for general $n$ in \cite{TWb}. 
 
 This paper is organized as follows. We construct moving frames and isotropic curvatures for isotropic curves in section \ref{ha}, give an explicit description of the tangent space of $\cm_{n+1, n}$ at $\g$ in section \ref{hf}, and review Drinfeld-Sokolov's construction of KdV type equations associated to the affine Kac-Moody algebras $\hat B_n^{(1)}$ and $\hat A_{2n}^{(2)}$ in section \ref{hb}. We prove results (a) and solve the Cauchy problem for these curve flows in section \ref{hc}. Result (b) is proved in section \ref{he}. We construct B\"acklund transformations for isotropic curve flows of B-type and A-type and give an algorithm to construct infinitely many families of explicit solutions of these isotropic curve flows in sections \ref{hd} and \ref{hg} respectively.  

\bs

\section{Moving frames along isotropic curves} \label{ha}

In this section, we prove the existence of isotropic parameter and construct isotropic moving frames and curvatures along isotropic curves. We also give some properties of the isotropic curvature map $\Psi$.  

The Lie algebra of $O(n+1, n)$ is 
\begin{align*}
o(n+1, n)&= \{A\in sl(2n+1, \R)\n A^t C_n + C_n A=0\}\\
&=\{(A_{ij}) \mid A_{ij}+(-1)^{i-j}A_{2n+2-j, 2n+2-i}=0, \quad 1 \leq i \leq 2n+1\}.
\end{align*}
Note that $A=(A_{ij})\in o(n+1, n)$ if and only if $A_{ij}$'s are symmetric (skew-symmetric resp.) with respect to the skew diagonal line $x+y= 2n+2$ if $i+j$ is odd (evenÄ resp.) and $A_{ij}=0$ if $i+j= 2n+2$. 
Let 
\beq\label{av}
 \cg_i=\sp\{e_{j, j+i} \mid 1 \leq i+j \leq 2n+1\}\cap o(n+1, n).
 \eeq
 Then
\begin{align*}
o(n+1, n)=\oplus_{i=-2n}^{2n}\cg_i, \quad \cg_{-2n}=\cg_{2n}={\bf 0}, \quad [\cg_i, \cg_j]\subset \cg_{i+j}.
\end{align*}

\bprop \label{da} \
\ben
\item[(i)] The dimension of a maximal isotropic subspace of $\R^{n+1, n}$ is $n$.
\item[(ii)] The $O(n+1, n)$-action on the space of ordered isotropic bases of $\R^{n+1, n}$ defined by $g\cdot (v_1, \ldots, v_n)= (gv_1, \ldots, gv_n)$ is transitive. 
\een
\eprop
\begin{proof} (i) Let $\{e_i, 1 \leq i \leq 2n+1\}$ denote the standard basis
of $\R^{2n+1}$. Then $A=\sp\{e_1, e_2, \cdots, e_n\}$ is an
isotropic space in $\R^{n+1,n}$. 

Let $V=\sp\{v_1, \cdots, v_n\}$ be another n-dimension isotropic subspace, $g_1=(e_1, \cdots, e_n)$, and $g_2=(v_1, \cdots, v_n)$. We claim that there exists $C \in O(n+1, n)$ such that $g_2=Cg_1$. From linear algebra, we can extend $\{v_1, \cdots, v_n\}$ to a basis $\{v_1, \cdots, v_n, v_{n+1}, \cdots, v_{2n+1}\}$ in $\R^{n+1, n}$ and denote $\ti{g}_2=(v_1, \cdots, v_{2n+1})\in O(n+1, n)$. Then choose $C=\ti{g}_2$.

(ii) Suppose $B=\sp\{w_1, \cdots, w_{n+1}\}$ is an isotropic subspace in $\R^{n+1, n}$ of dimension n+1. According to (i),  there exists $C \in O(n+1, n)$, such that $(w_1, \cdots, w_n)=C(e_1, \cdots, e_n)$. Therefore, we may assume $w_i=e_i, 1 \leq i \leq n$. Then from $\li e_i, w_{n+1} \ri=0$ for $1 \leq i \leq n$ and $\li w_{n+1}, w_{n+1} \ri=0$, we have $w_{n+1}=0$, which is a contradiction. This proves (ii). 
\end{proof}

\bprop
If $\gamma(s)$ is isotropic in $\R^{n+1, n}$ for all $s\in \R$, then there exists an orientation preserving parameter $x=x(s)$ unique up to translation such that $\li \g_x^{(n)}, \g_x^{(n)}
\ri = 1$, i.e.,  $x$  is the isotropic parameter of $\gamma$.
\eprop

\begin{proof} Since $\g$ is isotropic, 
$$\li \g_s^{(n-1)}, \g_s^{(i)} \ri=0, \quad 0 \leq i \leq n-1.$$ 
Take the derivative with respect to $s$ of both sides to get
$$
\li \g_s^{(n-1)}, \g_s^{(i)} \ri_s= \li \g_s^{(n)}, \g_s^{(i)}  \ri+\li \g_s^{(n-1)}, \g_s^{(i+1)} \ri=0
$$
So $\li \g_s^{(n)}, \g_s^{(i)}\ri=0$ for any $0 \leq i \leq n-1$. Since the span of $\{\g, \ldots, \g_s^{(n-1)}\}$ is maximal isotropic, $\li \g_s^{(n)}, \g_s^{(n)}\ri \not=0$. 

We claim that $\li \g_s^{(n)}, \g_s^{(n)}\ri >0$ for all $s\in \R$. 
To see this, we first note that Proposition \ref{da} (ii) implies that there exists $C \in O(n+1, n)$ such that 
$$C(\g, \cdots, \g_s^{(n-1)})=(e_1, e_2, \cdots, e_n),$$ 
where $e_i$ is the $i$-th standard basis of $\R^{2n+1}$. Let ${\bf c} =(c_1, c_2, \cdots, c_{2n+1})^t=C\g_s^{(n)}$. Then we have 
$$\li C\g_s^{(n)}, C\g_s^{(i)}\ri = \li \bc, e_i\ri = c_{2n+2-i}= \li \g_s^{(n)}, \g_s^{(i)}\ri =0$$ for $1\leq i\leq n$. So $c_{2n+2-i}=0$ for  $1 \leq i \leq n$. This implies that 
$$\li \g_x^{(n)}, \g_x^{(n)} \ri =\li C\g_x^{(n)}, C\g_x^{(n)} \ri=\bc^t C_n \bc= c_{n+1}^2.$$ 
But $\li \g_s^{(n)}, \g_s^{(n)}\ri \not=0$. This proves the claim. 

Choose $x$ such that
$\frac{d x}{d s}= \langle  \g_s^{(n)}, \g_s^{(n)} \rangle ^{1/2n}$ and the proposition follows. 
\end{proof}

\bthm\label{at}  
 Given $\gamma \in \cm_{n+1, n}$, there 
exists a unique isotropic moving frame $g=(p_1, \cdots, p_{2n+1})$ $:\R\to O(n+1, n)$ along $\g$, i.e.,  $p_i=\g_x^{(i-1)}$ for $1 \leq i \leq n+1$ and
$g^{-1}g_x=b+\sum_{i=1}^n u_i\b_i$ for some $n$ smooth functions $u_1, \cdots, u_n$, where $b$ and $\b_i$'s are given in \eqref{cj}. \ethm

\begin{proof}
We claim that $p_{n+2}, \cdots, p_{2n+1}$ and $u_1, \cdots, u_{n}$ can be constructed from recursive formulas. From $\li \g_x^{(n-1)}, \g_x^{(n)} \ri=0$, we have $\li \g_x^{(n-1)}, \g_x^{(n+1)} \ri=-\li \g_x^{(n)}, \g_x^{(n)} \ri=-1$. 

Let
\beq\label{ax}
u_1=-\frac{1}{2}\li \g_x^{(n+1)}, \g_x^{(n+1)} \ri,
\quad p_{n+2}=\g_x^{(n+1)}-u_1\g_x^{(n-1)}.
\eeq
Then 
$$
\bca
\li \g_x^{(i)}, p_{n+2} \ri=0, \quad 0 \leq i \leq n-2, \\
\li \g_x^{(n-1)}, p_{n+2}\ri=-1, \quad \li \g_x^{(n)}, p_{n+2} \ri=0, \quad \li p_{n+2}, p_{n+2} \ri=0.
\eca
$$
Since $\li p_{n+2}, \g_x^{(n)} \ri=0$, $\li (p_{n+2})_x, \g_x^{(n)} \ri+\li p_{n+2}, \g_x^{(n+1)} \ri =0$. From \eqref{ax},
$$\li (p_{n+2})_x, \g_x^{(n)}\ri= u_1. $$ 
On the other hand, $\li \g_x^{(n-2)}, p_{n+2} \ri=0$ implies that
$$\li \g_x^{(n-2)}, (p_{n+2})_x \ri=-\li \g_x^{(n-1)}, p_{n+2} \ri =1.$$
Set
\beq\label{ay}
u_2=\frac{1}{2}\li (p_{n+2})_x, (p_{n+2})_x \ri-\frac{1}{2}u_1^2,  \quad p_{n+3}=(p_{n+2})_x-u_2\g_x^{(n-2)}-u_1\g_x^{(n)}.
\eeq
Then 
$$
\bca
\li \g_x^{(i)}, p_{n+3} \ri =0, \quad 0 \leq i \leq n, i \neq n-2, \\
\li \g_x^{(n-2)}, p_{n+3}\ri=1, \quad  \li p_{n+2}, p_{n+3}\ri=0.
\eca
$$

Suppose we have already found $p_{n+2}, \cdots, p_{n+j}$ and $u_1, \cdots, u_{j-1}$ for $j\geq 3$ satisfying 
$$
\bca
\li p_i, p_{n+j} \ri =0, \quad 1 \leq i \leq n+j,  \quad i \neq n=2-j,\\
\li p_{n+2-j}, p_{n+j}\ri =(-1)^{j-1}, \\
(p_{n-1+j})_x=p_{n+j}+u_{j-1}p_{n-2+j}+u_{j-2}p_{n-4+j}.
\eca
$$
Set
$$
\bca
u_j=(-1)^{j}\li (p_{n+j})_x, (p_{n+j})_x \ri,  \\
p_{n+j+1}=(p_{n+j})_x-u_jp_{n+1-j}-u_{j-1}p_{n+3-j}.
\eca
$$
It follows from a direct computation that $\li p_{i}, p_{n+j+1} \ri=(-1)^{j}\d_{i, n+1-j}$. Hence we get $p_{n+j+1}$ and $u_j$.  The uniqueness of $g$ and $u_i$'s comes from the way they are constructed.
\end{proof}

The map $g$, $u_i$, and $u=\sum_{i=1}^n u_i \b_i$ in Theorem \ref{at} are the isotropic moving frame, the $i$-th isotropic curvature, and the isotropic curvature along $\g$. 

If follows from the Existence and Uniqueness Theorem of ordinary differential equations that we have 

\bprop \label{au} Let $\Psi:\cm_{n+1, n}\to V_n$ be the isotropic curvature map. Then $\Psi$ is onto and $\Psi^{-1}(u)$ is an $O(n+1, n)$-orbit for all $u\in C^\infty(\R, V_n)$. 
\eprop

\beg
Isotropic curves in $\R^{n+1, n}$ with zero isotropic curvatures are of the form
\beq
\gamma=C(1,x, \frac{x^2}{2!}, \cdots, \frac{x^{2n-1}}{(2n-1)!},\frac{x^{2n}}{2n!})^t, \ \ \ C\in O(n+1, n).
\eeq
\eeg

\beg\label{dt}
The isotropic moving frame $g=(\g, \g_x, p_3)$ and isotropic curvature $q:=u_1$ along  $\g \in \cm_{2,1}$ are 
\begin{align*}
&q=-\frac{1}{2}\li\g_{xx},\g_{xx} \ri,\quad  p_3=\g_{xx}-q\g,\\
& g^{-1}g_x=\left(\begin{array}{ccc}0&q&0
\\1&0&q\\0&1&0\end{array}\right).
\end{align*}
\eeg

\beg\label{gb}
For $\g \in \cm_{3,2}$, let $u_1=-\frac{1}{2} \li \g_{x}^{(3)}, \g_{x}^{(3)}\ri$, $p_4=\g_x^{(3)}-u_1\g_x$, $u_2=\frac{1}{2}(\li (p_4)_x, (p_4)_x \ri-u_1^2)$, and $p_5=(p_4)_x-u_2\g-u_1\g_{xx}$. Then $g=(\g, \g_x, \g_{xx}, p_4, p_5)$ is the isotropic moving frame along $\g$, i.e.,  
$$g^{-1}g_x=\bpm 0 & 0 & 0& u_2 & 0 \\ 1 & 0 & u_1 & 0 & u_2  \\
0 & 1 & 0& u_1& 0 \\ 0 & 0& 1 & 0 & 0 \\ 0 & 0 & 0& 1 & 0\epm,$$
and $u_1, u_2$ are the isotropic curvatures of $\g$. 
\eeg

\bprop \label{ao}
Let $\Psi: \cm_{n+1, n} \rightarrow C^{\infty}(\R, V_n)$ be the isotropic curvature map. Then the differential of $\Psi$ at $\g$ is 
\beq\label{ce}
d\Psi(\d \g)=\d u =[\p_x+b+u, g^{-1}\d g],
\eeq
where $g$, $u$, and $\d g$ are the isotropic moving frame, isotropic curvature, and the variation of $g$ when we vary $\g$ by $\d \g$ respectively.
\eprop
\begin{proof} It follows from $g^{-1}g_x=b+u$ that we have 
$$\d u=-g^{-1}\d g g^{-1}g_x+g^{-1}(\d g)_x=-g^{-1}\d g(b+u)+g^{-1}(\d g)_x,$$
On the other hand,
$$(g^{-1}\d g)_x=-g^{-1}g_xg^{-1} \d g+g^{-1}(\d g)_x=-(b+u)g^{-1}\d g+g^{-1}(\d g)_x,$$
Therefore
$$\d u=-g^{-1}\d g(b+u)+(g^{-1}\d g)_x+(b+u)g^{-1}(\d g)=[\p_x+b+u, g^{-1}\d g].$$
\end{proof}

As a consequence of Propositions \ref{au} and \ref{ao} we have

\bcor \label{es}
Given $u,\eta\in C^\infty(\R, V_n)$, then there exists $\xi\in C^\infty(\R, o(n+1, n))$ such that $[\p_x+ b+u, \xi]= \eta$, where $b$ is as in \eqref{cj}.
\ecor

\bs
\section{The tangent space of $\cm_{n+1, n}$ at $\g$} \label{hf}

In this section, we identify the tangent space of $\cm_{n+1, n}$ at $\g$ as $C^\infty(\R, \R^n)$.

Henceforth in this paper we set
$$e_1=(1, 0, \ldots, 0)\in \R^{n+1, n}.$$
Below we give a useful description of $T_\g\cm_{n+1, n}$. 

\bprop\label{bq} Let $g$ and $u$ be the isotropic moving frame and isotropic curvature along $\g\in \cm_{n+1,n}$. Suppose $C:\R\to O(n+1, n)$ satisfies 
\beq\label{ae}
[\p_x+b+u, C] \in C^{\infty}(\R, V_n).
\eeq 
 Then $\xi(\g)=gCe_1$ is tangent to $\cm_{n+1, n}$ at $\g$. Conversely, all tangent vectors of $\cm_{n+1, n}$ at $\g$ arise this way. 
  \eprop
  
\begin{proof} It follows from the definition of $\cm_{n+1,n}$ that $\d \g$ is tangent to $\cm_{n+1, n}$ at $\g$ if and only if 
\beq\label{cn}
\bca
\li (\d \g)_x^{(i)}, \g_x^{(j)} \ri+\li \g_x^{(i)}, (\d \g)_x^{(j)} \ri=0, \quad 0 \leq i, j \leq n-1, \\
\li (\d \g)_x^{(n+1)}, \g_x^{(n)} \ri=0. 
\eca
\eeq
Let $\eta_j$ denote the $j$-th column of $gC$ for $1\leq j\leq 2n+1$. To prove $gCe_1$ is tangent to $\cm_{n+1, n}$, it suffices to prove that $\eta_1$ satisfies \eqref{cn}.
Let $\rho=[\p_x+b+u, C]$. A direct computation gives 
$$(gC)_x=g_xC+gC_x=gC(b+u)+g\rho.$$
Since the first $n$ columns of $\rho$ are zero, the first $n+1$ columns of $gC$ are related by
$$\eta_2=(\eta_1)_x, \cdots, \eta_{n+1}=(\eta_1)_x^{(n)}.$$
Hence, for $0 \leq i, j \leq n-1$, we have 
\begin{align*}
 & \li (\eta_1)_x^{(i)}, \g_x^{(j)} \ri+\li \g_x^{(i)}, (\eta_1)_x^{(j)} \ri = \li gCe_{i+1}, ge_{j+1} \ri+\li ge_{i+1}, gCe_{j+1} \ri \\
& \quad  = \li Ce_{i+1}, e_{j+1} \ri +\li e_{i+1}, Ce_{j+1} \ri  = e_{i+1}^t(C^tC_n+C_nC^t)e_{j+1} \\
& \quad =0.
\end{align*}
Since $C=(C_{ij}) \in o(n+1, n)$, $\li (\eta_1)_x^{(n)}, \g_x^{(n)} \ri=C_{n+1, n}=0$. So $\xi (\g)=\eta_1$ is tangent to $\cm_{n+1, n}$ at $\g$.

The converse follows from Proposition \ref{ao}. 
\end{proof}

Next we prove that if $C=(C_{ij})$ satisfying \eqref{ae} then $C$ is determined by $\{C_{n+i, n+1-i}, 1 \leq i \leq n\}$ or $\{C_{2i, 1}, 1 \leq i \leq n\}$.

\bthm \label{ag} 
Let $u \in C^{\infty}(\R, V_n)$, $C=(C_{ij}) \in C^{\infty}(\R, o(n+1, n))$ and $v_i=C_{n+i, n+1-i}$ for $ 1 \leq i \leq n$. Suppose $[\p_x+b+u, C] \in V_n$. Then we have the following.
\ben
\item[(i)] $C_{ij}$'s are differential polynomials in $u, v_1, \cdots, v_n$.
\item[(ii)] $C_{2n-2i, 1}=v_{n-i}+\phi_i$ for $0 \leq i \leq n-1$, where $\phi_i$ is a differential polynomial in $u, v_{n+1-i}, \cdots, v_n$.
\item[(iii)] There exist differential polynomials $h_{2i+1}$ such that 
$$C_{2i+1, 1}=h_{2i+1}(u, C_{2i+2, 1}, \cdots, C_{2n, 1}), \quad 0 \leq i \leq n-1.$$  
\item[(iv)] $C_{ij}$'s are differential polynomials of $u, C_{21}, C_{41}, \cdots, C_{2k, 1}, \cdots, C_{2n,1}.$
\een
\ethm

\begin{proof}  Since $C \in o(n+1, n)$, $C_{n+1+i, n+2-i}=C_{n+i, n+1-i}=v_i$. Let  $v=\sum_{i=1}^nv_i\b_i^t \in V_n^t$, where $V_n$ and $\b_i$'s are defined in \eqref{aj}. Let $\cg_k$ be as in \eqref{av}. Then $[\cg_i, \cg_j] \subset \cg_{i+j}$. For $\xi \in o(n+1, n)$, let $\xi_{\cg_i}$ denote the $\cg_i$-component of $\xi$ w.r.t. $o(n+1, n)=\oplus_{i=1-2n}^{2n-1}\cg_i$. 

Suppose $[\p_x+b+u, C]=\sum_{i=1}^{n} \eta_i\b_i$. Write $C=\sum_{i=1-2n}^{2n-1}C_i$ with $C_i \in \cg_i$. Then 
\beq\label{by}
C_j'+[b, C_{j+1}]+[u, C]_{\cg_j}=\bca
\eta_{i}\b_i,  \quad j=2i-1, \\
0, \quad \text{else}.
\eca
\eeq
We claim that $C_j$ are differential polynomials in $v$ and $u$. For $j=1-2n$, we have $C_{2n, 1}=C_{2n+1, 2}=v_{n}$. For $j < 0$, if $j$ is even, $ad(b): \cg_j \rightarrow \cg_{j-1}$ is a bijection. If $j$ is odd, then $\dim(\Im(ad(b)(\cg_j)))=\dim(\cg_{j-1})=\dim(\cg_j)-1$. Then from \eqref{by}, for both cases, entries of $C_j$ are differential polynomials in $v_n, \cdots, v_{-[\frac{j}{2}]}$. Then by induction, the claim is true for $j < 0$.

Note that $ad(b)$ is a bijection from $\cg_0$ to $\cg_{-1}$, and we have the $\cg_j$ component $[u, C]_{\cg_j}$ depends only on $u$, $v_1, \cdots, v_n$. So $C_0$ is a differential polynomial in $u$ and $v$. 

For $j>0$, we see that when $j$ is odd, $ad(b): \cg_j \rightarrow \cg_{j-1}$ is again a bijection. When $j$ is even, we have $\dim(\Im(ad(b)(\cg_j)))=\dim(\cg_{j})=\dim(\cg_{j-1})-1$. Therefore, in both cases, $C_j$ can be solved uniquely from $C_{j-1}$ and $\eta_i$'s are differential polynomials in entries of $C_{2i-1}$. By induction,  the claim is true for $j > 0$. This proves the statement (i).
 
To prove (ii), let $j=2i+1-2n$ in \eqref{by}. Then the linear system implies that $C_{2n-2i, 1}=v_{n-i}+\phi_i, 0 \leq i \leq n-1$, where $\phi_i$ is a differential polynomial in $u, v_{n+1-i}, \cdots, v_n$.

Statement (iii) and (iv) are consequence from (i) and (ii).
\end{proof}

\bcor \label{az} Let  $g=(p_1, \ldots, p_{2n+1})$ and $u$ denote the isotropic moving frame and isotropic curvature along $\g \in \cm_{n+1, n}$. Then $\xi=\sum_{i=1}^{2n}\xi_i p_i$ is tangent to $\cm_{n+1, n}$ at $\g$ if and only if 
\beq\label{ee}
\xi_{2i+1}=h_{2i+1}(u, \xi_{2i+2}, \cdots, \xi_{2n}), \quad 0 \leq i \leq n-1,
\eeq
where $h_{2i+1}$'s are the differential polynomials given in Theorem \ref{ag}. In particular, we identify $T_\g\cm_{n+1, n}$ as $C^\infty(\R, \R^n)$. 
\ecor

The proof of Theorem \ref{ag} gives the following.

\bcor \label{cz} Let $u \in C^{\infty}(\R, V_n)$, and $v_1, \cdots, v_n \in C^{\infty}(\R, \R)$. Then there exists  a unique $C=(C_{ij}) \in C^{\infty}(\R, o(n+1, n))$ such that $[\p_x+b+u, C] \in V_n$ and $C_{n+i, n+1-i}=v_i$ for $1 \leq i \leq n$. 
\ecor

\bdefn \label{aw}
Given  $u \in C^{\infty}(\R, V_n)$, let
$$P_u: C^{\infty}(\R, V_n^t) \rightarrow C^{\infty}(\R, o(n+1, n))$$
be the map defined by $P_u(v)=C$ for $v=\sum_{i=1}^n v_i \b_i^t$, where $C$ is the unique $o(n+1, n)$-value map defined by $v_1, \ldots, v_n$ in Corollary \ref{cz}.
\edefn

\bcor\label{eh} Let $g$ and $u$ be the isotropic moving frame and the isotropic curvature  along  $\g \in \cm_{n+1, n}$. Then the following statements are equivalent for $C:\R\to o(n+1, n)$:
\ben 
\item $gCe_1$ is tangent to $\cm_{n+1, n}$ at $\g$.
\item $[\p_x+b+u, C]\in V_n$.
\item $C= g^{-1}\d g$, where $\d g$ is the  variation of $g$ when we vary $\g$ by $\d \g$.
\item $C=P_u(\pi_0(C))$,
\een
where $\pi_0: o(n+1, n)\to V_n^t$ be the projection defined by
\beq\label{bk}
\pi_0((y_{ij}))=\sum_{i=1}^{n}y_{n+i, n+1-i}\b_i^t, \quad \b_i^t=e_{n+i, n+1-i}+e_{n+1+i, n+2-i}.
\eeq
\ecor

\bcor\label{gk} Given $\g\in \cm_{n+1,n}$ and $\xi=(\xi_1, \ldots, \xi_{2n},0)^t$ satisfying \eqref{ee}, then there is a unique $o(n+1, n)$-valued differential polynomial $\Theta_\g(\xi)\in C^\infty(\R, o(n+1, n))$ such that the $i1$-th entry is $\xi_i$ for $1\leq i\leq 2n+1$ and $[\p_x+b+u, \Theta_\g(\xi)]\in C^\infty(\R, V_n)$. Moreover, $\Theta_\g(\xi)= P_u(\pi_0(\Theta_\g(\xi)))$. 
\ecor

Below we give some examples of $T(\cm_{n+1, n})_\g$ when $n$ is small.

\beg \label{dk} {\bf $T_\g\cm_{2, 1}$} \ 

When $n=1$, $b=e_{21}+e_{32}$, $\b_1=e_{12}+e_{23}$, $V_1=\R \b_1$, $u=q\b_1$. a direct computation implies that $[\p_x+b+u, C]\in V_1$ if and only if 
$$C=\bpm -\xi_x & -\xi_{xx}+q\xi & 0  \\
\xi & 0 & -\xi_{xx}+q\xi \\ 0 & \xi & \xi_x \epm.$$
Moreover,  $[\p_x+b+u, C]=(-\xi_{x}^{(3)}+2q\xi_x+(q)_x\xi)\b_1$ and  
$$T(\cm_{2, 1})_\g=\{-\xi_x\g+\xi\g_x \mid \xi \in C^{\infty}(\R, \R)\}.$$
\eeg

\beg\label{dl} {\bf $T_\g\cm_{3, 2}$} \ 
 
Here $b=\sum_{i=1}^4e_{i+1, i}$, $\b_1=e_{23}+e_{34}$, $\b_2=e_{14}+e_{25}$, $V_2=\R\b_1\oplus \R \b_2$, and $u=u_1\b_1+u_2\b_2$. Then $C \in C^{\infty}(\R, o(3, 2))$ satisfying $[\p_x+b+u, C]\in V_2$ if and only if 
$$C=\bpm \eta_x^{(3)}-2\xi_x-(u_1)_x\eta& \ast & \ast & \ast & 0  \\ \xi & a & \ast & 0 & \ast \\ -(\eta)_x & \zeta & 0 & \ast & \ast \\ \eta & 0 & \zeta  & -a & \ast \\ 0& \eta & (\eta)_x & \xi & -\eta_x^{(3)}+2\xi_x+(u_1)_x\eta\epm$$
for some real valued functions $\xi, \eta$,
where $a=\eta_x^{(3)}-\xi_x-(u_1\eta)_x$ and $\zeta=\xi-\eta_{xx}+u_1\eta$.  So
$$T(\cm_{3, 2})_\g=\{(\eta_x^{(3)}-2\xi_x-(u_1)_x\eta)\g+\xi\g_x-\eta_x\g_{xx}+\eta p_4 \mid \xi, \eta \in C^{\infty}(\R, \R)\},$$
where $g=(\g, \g_x, \g_{xx}, p_4, p_5) \in O(3 ,2)$ is the isotropic moving frame along $\g$.

In particular, we see that $X_1(\g)=-\frac{1}{2}u_1\g_x+p_4$ and $X_2(\g)=-\frac{3}{5}(u_1)_x\g-\frac{1}{5}u_1\g_x+p_4$ are tangent vector fields of $\cm_{3,2}$. Hence
\beq\label{ap}
\g_t=-\frac{1}{2}u_1\g_x+p_4,
\eeq
\beq\label{aq}
\g_t=-\frac{3}{5}(u_1)_x\g-\frac{1}{5}u_1\g_x+p_4,
\eeq
are isotropic curve flows on $\cm_{3,2}$. By Example \ref{gb}, we have  $p_4= \g_x^{(3)}- u_1 \g_x$. So \eqref{ap} and \eqref{aq} are \eqref{ar} and \eqref{as} with $n=2$ respectively.  
\eeg

\beg {\bf $\cm_{n+1, n}$, ($n \geq 3$)} \ 

For $\g \in \cm_{n+1, n}$, a direct computation implies that 
$$X(\g)=-(n\xi_x+3(u_1)_x)\g+\xi\g_x+\g_x^{(3)}$$ is tangent to $\cm_{n+1, n}$ at $\g$. So 
\beq\label{af}
\g_t=-(n\xi_x+3(u_1)_x)\g+\xi\g_x+\g_x^{(3)}
\eeq 
is an isotropic curve flow on $\cm_{n+1, n}$. 
For example, if we choose $\xi$ to be $-\frac{3}{n} u_1$ and $ -\frac{6}{2n+1} u_1$ then we get \eqref{ar} and \eqref{as} given in the introduction respectively. 
\eeg

\bs
\section{The $\bn1$-KdV and $\an2$-KdV hierarchies}\label{hb}

Drinfeld-Sokolov associated to each affine Kac-Moody algebra a KdV-type soliton hierarchy  (cf. \cite{DS84}). These hierarchies are constructed as quotient flows for some gauge group actions. Different cross sections of the gauge group action give different but equivalent hierarchies. In this paper, we construct a suitable cross section for the gauge action so that the differential invariants for the isotropic curves lies in this cross section.  We also prove that the $\an2$-KdV hierarchy is the constraint KP hierarchy in \cite{YC92}. 

Let $\C^{n+1, n}$ be the vector space of $\C^{2n+1}$ equipped with the bilinear form $\li X, Y \ri= X^tC_nY$, where $C_n$ is as in \eqref{am}. Let $O_{\C}(n+1, n)$ be the Lie group preserving
$\langle \ , \ \rangle$ on $\C^{n+1, n}$, i.e.
\beq\label{bw}
O_{\C}(n+1, n)=\{A \in SL(2n+1, \C) \mid  A^tC_nA=C_n\}.
\eeq
Its Lie algebra is  $o_{\C}(n+1, n)=\{A \in sl(2n+1, \C) \mid  A^tC_n+C_nA=0\}$.
Let $\cb_n^+, \cn_n^+, \ct_n$ denote the subalgebras of upper triangular, strictly upper triangular, and diagonal matrices in $o(n+1, n)$ respectively. Let $B_n^+$ and $N_n^{+}$ be the connected subgroups of $O(n+1, n)$ with Lie algebras $\cb_n^+$ and $\cn_n^+$.  

\ms
First we review the construction of the {\bf $\bn1$-KdV hierarchy\/}. Set  
\begin{align*}
& \bn 1 =\left\{\xi(\l)=\sum_{i \leq m_0}\xi_i \l^i  \n \xi_i \in o_\C(n+1, n), \overline{\xi(\bar{\l})}=\xi(\l), m_0\in \Z\right\}, \\
& (\bn 1)_+ = \left\{\xi(\l)=\sum_{i \geq 0} \xi_i \l^i \in \bn 1\right\}, (\bn 1)_-=\left\{\xi(\l)=\sum_{i < 0}\xi_i\l^i \in \bn 1\right\}. \notag
\end{align*}
Note that 
\ben
\item[(i)] $\xi(\l)=  \sum_i \xi_i \l^i\in \bn1$ if and only if $\xi_i \in o(n+1, n)$ for all $i$,
\item[(ii)] $\bn 1=(\bn 1)_+\oplus (\bn 1)_-$ is a direct sum of linear subspaces.
\een

Let 
\beq\label{db}
\b=\frac{1}{2}(e_{1, 2n}+e_{2, 2n+1})
\eeq
and 
\beq\label{dc}
J_B=\sum_{i=1}^{2n}e_{i+1, i}+\frac{1}{2}\l(e_{1, 2n}+e_{2, 2n+1})=b+\beta \l.
\eeq
Note that $J_B^{2j}\not\in \bn1$ and 
\beq\label{dh}
J_{B}^{2j-1} (j \geq 1) \in (\bn 1)_+, \quad J_B^{2n+1}=\l J_B.
\eeq
So we can use the splitting $\bn1= (\bn1)_+\oplus (\bn1)_-$ and the commuting sequence $\{J_B^{2j-1}\n j\geq 1\}$ in $(\bn1)_+$ to construct a hierarchy of soliton equations (cf. \cite{TU11}).  A direct computation implies that given $q \in C^{\infty}(\R, \cb_n^+)$, there
exists a unique $P(q, \l)$ conjugate to $J_B$ and satisfying
\beq\label{ch}
[\p_x+J_B+q, P(q, \l)]=0,
\eeq
Expand $P^{2j-1}(q,\l)$ as a power series in $\l$, 
\beq\label{dg}
P^{2j-1}(q, \l)=\sum_{i \leq [\frac{2j-1}{2n+1}]+1}P_{2j-1, i}(q)\l^i.
\eeq
It can be checked that all $P_{2j-1, i}(q)$'s are differential polynomials of $q$.

The $(2j-1)$-th flow ($j \geq 1$) on $C^{\infty}(\R, \cb_n^+)$  is 
\beq\label{ah}
q_{t_{2j-1}}=[\p_x+b+q, P_{2j-1, 0}(q)].
\eeq
These flows commute. Moreover, $q$ is a solution of \eqref{ah} if and only if
$$[\p_x+J_B+q, \p_{t_{2j-1}}+(P^{2j-1}(q, \l))_+]=0$$
for any $\l \in \C$. Here $\xi_+$ is the projection of $\bn 1$ onto $(\bn 1)_+$ along $(\bn 1)_-$. 

\bdefn
We call $F(x, t, \l)\in O_\C(n+1, n)$ a \emph{frame} of the solution $q$ of \eqref{ah} if $F$ is holomorphic for $\l\in \C$ and satisfies
$$F^{-1}F_x=J_B+q, \quad F^{-1}F_t=(P^{2j-1}(q, \l))_+, \quad \overline{F(x, t, \bar \l)}=F(x, t, \l).$$
\edefn

The group $C^{\infty}(\R, N_n^+)$ acts on $C^{\infty}(\R, \cb_n^+)$ by gauge transformations, 
\beq\label{dq}
\p_x+J_B+ \D \ast q= \D(\p_x+J_B+ q)\D^{-1},
\eeq or equivalently,
\beq\label{dr}
\D \ast q=\D (J_B+q)\D^{-1}-\D_x\D^{-1}-J_B \in \cb_n^+.
\eeq

A direct computation similar to the one given in \cite{TU11} implies the following.

\bprop \

\ben
\item[(i)] Let $q\in C^\infty(\R, \cb_n^+)$, and $\D\in C^\infty(\R, N_n^+)$. Then 
$$\D P(q,\l) \D^{-1} = P(\D\ast q, \l),$$
where $P(q,\l)$ is defined by \eqref{ch}. 
\item[(ii)] Let $q$ be a solution of \eqref{ah}, $\D\in C^\infty(\R, N_n^+)$, and $\ti q(\cdot, t)= \D\ast q(\cdot, t)$. Then $\ti q$ is again a solution of \eqref{ah}. 
\een
\eprop

So \eqref{ah} is invariant under the action of $C^\infty(\R, N_n^+)$ and it induces a quotient flow on the orbit space $\frac{C^\infty(\R, \cb_n^+)}{C^\infty(\R, N_n^+)}$.  The next Proposition, which can be proved by a direct computation, shows that $C^\infty(\R, V_n)$ is a cross section of this gauge action. 

    \bprop \label{ab}
Given $q \in C^{\infty}(\R, \cb_n^+)$, then there
exist a unique $\D \in C^{\infty}(\R, N_n^+)$ and $u\in C^\infty(\R, V_n)$ such that
\beq\label{aa}
\D(\p_x+J_B+q)\D^{-1}=\p_x+J_B+u, 
\eeq
where $V_n$ is as in \eqref{aj} and $J_B$ is given by \eqref{dc}.  Moreover, entries of $\D$ and $u$ are differential polynomials of $q$.
\eprop

\bdefn\label{aba} Let $\G: C^{\infty}(\R, \cb_n^+) \rightarrow C^{\infty}(\R, V_n)$ and $D: C^{\infty}(\R, \cb_n^+) \rightarrow C^{\infty}(\R, N_n^+)$ be the maps defined by $\G(q)=u$ and $D(q)=\D$, where $q, u$, and $\D$ are related by \eqref{aa} as in Proposition \ref{ab}.
\edefn

\bcor \label{cg} (\cite{DS84}) Let $u\in C^\infty(\R, V_n)$, $j \geq 1$, and $P_{2j-1,0}(u)$ be defined by \eqref{dg}. 
 Then there exists a unique
differential polynomial $\eta_j \in C^{\infty}(\R, \cn_n^+)$ of $u\in C^\infty(\R, V_n)$ such that
\beq\label{bp}
[\p_x+b+u, P_{2j-1,0}(u)-\eta_j(u)]\, \in C^{\infty}(\R, V_n).
\end{equation}
\ecor
 
   The $(2j-1)$-th flow \eqref{ah} induces a quotient flow on the cross section $C^\infty(\R, V_n)$ by projecting the solutions to the cross section along orbits. This quotient flow is the  \emph{ $(2j-1)$-th $\bn 1$-KdV flow\/},
   \beq\label{ai}
u_{t_{2j-1}}= [\p_x+b+u, P_{2j-1,0}(u)-\eta_j(u)]
\eeq
 
\bcor\label{cr} Let $\pi_0:o(n+1, n)\to V_n^t$ be the projection defined by \eqref{bk}, and $P_u$ the operator defined in Definition \ref{aw}. Then
\ben
\item[(i)] $P_{2j-1,0}(u)-\eta_j(u)= P_u(\pi_0(P_{2j-1,0}(u))$, $\pi_0(P_{2j-1,0}(u) - \eta_j(u)) $ $=\pi_0(P_{2j-1,0}(u))$, 
\item[(ii)] the $(2j-1)$-th $\bn1$-KdV flow \eqref{ai} can be written as 
\beq\label{aia} 
u_{t_{2j-1}}= [\p_x+b+u, P_u(\pi_0(P_{2j-1,0}(u)))], 
\eeq
\een
where $P_{2j-1, 0}(u)$ is given in \eqref{dg}. 
\ecor

\begin{proof}
 Since  $\eta_j(u)$ is strictly upper triangular,  we have 
  $$\pi_0(P_{2j-1,0}(u)-\eta_j(u))= \pi_0(P_{2j-1, 0}(u)).$$ 
 Then (i) and (ii) follows from \eqref{bp} and Corollary \ref{eh}. 
\end{proof}

It follows from the construction of the $\bn1$-KdV flows that we have the following.

\bprop\label{dv}
 Let $u$, $P(u,\l)$, $P_{2j-1,0}(u)$, and $\eta_j(u)$ be the differential polynomials given in Corollary \ref{cg}. Then the following statements are equivalent:
\ben
\item[(i)] $u$ is a solution of the $(2j-1)$-th $\bn1$-KdV flow \eqref{ai}.
\item[(ii)] $[\p_x+ b+u, \p_{t_{2j-1}}+ P_{2j-1,0}(u)-\eta_j(u)]=0$.
\item[(iii)]
$[\p_x+J_B+u, \, \p_{t_{2j-1}}+ (P^{2j-1}(u,\l))_+- \eta_j(u)]=0$ for all parameters $\l \in \C$ (this is the {\it Lax pair} for \eqref{ai}).
 \item[(iv)] The linear system 
 \beq\label{eo}
 \bca g^{-1}g_x= b+u,\\ g^{-1}g_t= P_{2j-1,0}(u)-\eta_j(u),\eca
 \eeq is solvable for $g:\R^2\to O(n+1, n)$.
 \item[(v)] For $\l\in \C$, the linear system 
 \beq\label{en}
 \bca E^{-1}E_x= J_B+ u, \\ E^{-1}E_t= (P^{2j-1}(u,\l))_+- \eta_j(u),\eca
 \eeq is solvable for $E(\cdot,\cdot, \l): \R^2 \rightarrow O_\C(n+1, n)$.
 \een
 \eprop

\bdefn We call $E(x,t,\l)\in O_\C(n+1, n)$ a {\it frame\/} of the solution $u$ of \eqref{ai} if $E(x,t,\l)$ is holomorphic for $\l\in \C$ and is a solution of \eqref{en} satisfying 
$$\overline{(E(x,t,\bar \l))}=E(x, t, \l).$$
\edefn

The next two Propositions follow from the constructions of flows \eqref{ah} and \eqref{ai}. 

\bprop\label{cga} Let $q: \R^2 \rightarrow \cb_n^+$ be a solution of \eqref{ah}, and $\D(\cdot, t)=D(q(\cdot, t))$, where $D$ is the operator defined in Definition \ref{aba}. Then $u=\D\ast q$ is a solution of the $(2j-1)$-th $\bn 1$-KdV flow \eqref{ai}, where the action $\ast$ is defined by \eqref{dr}.   Moreover, if $F$ is a frame of the solution $q$ of \eqref{ah}, then $E=F \D^{-1}$ is a frame of the solution $u$ of \eqref{ai}.
\eprop

\bprop\label{do} Let $u$ be a solution of \eqref{ai}. Suppose $\D:\R^2\to N_n^+$ satisfying $\D_t\D^{-1}=\eta_j(u)$. Then $q(\cdot, t)= \D(\cdot, t)^{-1}\ast u(\cdot, t)$ is a solution of \eqref{ah}, where the action $\ast$ is defined by \eqref{dr}.  Moreover, if $E$ is a frame of the solution $u$ of \eqref{ai}, then $F= E\D$ is a frame for the solution $q$ of \eqref{ah}.
\eprop

Note that if $\D_t(x, t)\D^{-1}(x, t)= \eta_j(u)$ and $f(x) \in C^\infty(\R, N_n^+)$, then $\ti \D(x, t)= \D (x, t) f(x)$ also satisfies $\ti \D_t \ti \D^{-1}=\eta_j(u)$ and $\ti q= \ti \D^{-1}\ast u= f\ast q$ is again a solution of \eqref{ah}. 

Next we use \eqref{ch} and \eqref{bp} to write down the $(2j-1)$-th $\bn1$-KdV flows for small $n$ and $j$.

\beg\label{ea} {\bf The $\hat B_1^{(1)}$-KdV hierarchy} \

In this case, $J_B=\left(\begin{array}{ccc}0&\frac{\l}{2}&0\\1&0&\frac{\l}{2}\\0&1&0\end{array}\right)$ and  $u=\left(\begin{array}{ccc}0&q&0\\0&0&q\\0&0&0\end{array}\right)$.  A direct computation implies that 
\begin{align*}
&P_{3, 0}(u)= \bpm q_x & 0 & 0 \\ -q & 0 & 0 \\ 0 & -q & -q_x \epm,\quad \eta_{2}(u)=\bpm 0 & q^2-q_{xx} & 0 \\ 0 & 0  & q^2-q_{xx}  \\ 0 & 0 & 0 \epm, \\
& (P^3(u, \l))_+-\eta_2(u)=J_B^3+\left(\begin{array}{ccc}q_x&q_{xx}-q^2&0\\-q&0&q_{xx}-q^2\\0&-q&-q_x\end{array}\right).
\end{align*}
The third $\hat B_{1}^{(1)}$-KdV flow is the KdV $q_t=q_{xxx}-3qq_x$. Since $sl(2, \R)$ is isomorphic to $o(2, 1)$, the algebra $\hat B_1^{(1)}$ is isomorphic to $\hat A_{1}^{(1)}$. So the $\hat B_1^{(1)}$-KdV hierarchy is the KdV hierarchy under this isomorphism. 
\eeg

\beg {\bf The $\hat B_2^{(1)}$-KdV hierarchy} \

We have $J_B=\sum_{i=1}^{4}e_{i+1, i}+\frac{\l}{2}(e_{14}+e_{25})$ and
$u=u_1(e_{23}+e_{34})+u_2(e_{14}+e_{25})$.
Then
\begin{align*}
& P_{3,0}(u) =  \\
&\bpm 0& u_2 & u_2' & u_2''-\frac{1}{2}u_1u_2 & 0 \\
-\frac{1}{2}u_1 & -\frac{1}{2}u_1' & \frac{1}{2}(u_1^2-u_1'')+2u_2& 0 & u_2''-\frac{1}{2}u_1u_2 \\
0& \frac{1}{2}u_1 & 0 & \frac{1}{2}(u_1^2-u_1'')+2u_2 & -u_2' \\
1& 0& \frac{1}{2}u_1 & \frac{1}{2}u_1' & u_2 \\
0& 1 & 0& -\frac{1}{2}u_1& 0\epm,
\end{align*} 
and 
\begin{align*}
&  P_{3,0}(u)-\eta_2(u) =  \\
& \bpm 0& u_2 & u_2' & u_2''-\frac{1}{2}u_1u_2 & 0 \\
-\frac{1}{2}u_1 & -\frac{1}{2}u_1' & \frac{1}{2}(u_1^2-u_1'')+2u_2& 0 & u_2''-\frac{1}{2}u_1u_2 \\
0& \frac{1}{2}u_1 & 0 & \frac{1}{2}(u_1^2-u_1'')+2u_2 & -u_2' \\
1& 0& \frac{1}{2}u_1 & \frac{1}{2}u_1' & u_2 \\
0& 1 & 0& -\frac{1}{2}u_1& 0\epm.
\end{align*}
So the third flow is
\begin{equation*}
\begin{cases}
(u_1)_t=-\frac{1}{2}u_1^{(3)}+\frac{3}{2}u_1(u_1)_x+3(u_2)_x,\\
(u_2)_t=u_2^{(3)}-\frac{3}{2}u_1(u_2)_x.
\end{cases}
\end{equation*}
\eeg

\beg\label{el}  Although we do not have the explicit formula for $P_{3,0}(u)$ for general $n\geq 3$, a direct computation implies that the first column of $P_{3,0}(u)$ is $(0, -\frac{3}{n} u_1, 0, 1, 0,  \ldots, 0)^t$. 
\eeg

\ms
Next we construct the {\bf  $\an2$-KdV hierarchy}.  
Let 
\begin{align*}
& \anone=\left\{\xi(\l)=\sum_{i \leq n_0}\xi_i\l^i \n \xi_i \in sl(2n+1, \C), \overline{\xi (\bar \l)}=\xi(\l) \right\},\\
& (\anone)_+=\left\{\sum_{i \geq 0}\xi_i\l^i \in \anone\right\}, \quad (\anone)_-=\left\{\sum_{i < 0}\xi_i\l^i \in \anone\right\}.
\end{align*}
Let $C_n$ be as in \eqref{am}. Set  
\begin{align*}
& \an2 =\{\xi(\l) \in \anone \n C_n \xi^t(-\l)C_n+\xi(\l) =0\},\\
& (\an2)_{\pm}= \an2 \cap (\anone)_\pm.
\end{align*}
Then $\an2=(\an2)_+ \oplus (\an2)_-$. Note that $\xi(\l)=\sum_{i \leq n_0}\xi_i\l^i \in \an 2$ if and only if $\xi_{2i} \in o_{\C}(n+1, n)$ and $C_n\xi_{2i+1}^tC_n=\xi_{2i+1}$.

Set
\beq\label{df}
J=\left(\sum_{i=1}^{2n}e_{i+1, i}\right)+e_{1, 2n+1}\l=b+e_{1, 2n+1}\l.
\eeq
Then $J^{2n+1}=\l \I_{2n+1}$, $J^{2i}\not\in \an2$, and $J^{2j-1}\in (\an2)_+$. So the splitting of $\an2$ and $J^{2j-1}$ produce $(2j-1)$-th flow with $(2j-1)\not\equiv 0$ $({\rm mod} (2n+1))$.  

A direct computation implies that given $q \in C^{\infty}(\R, \cb_n^+)$, there
exists a unique $S(q, \l) \in \hat A_{2n}^{(2)}$ such that
\beq\label{aka}
\bca
[\p_x+J+q, S(q, \l)]=0, \\
S(q, \l) \text{\ is conjugate to\ } J.
\eca
\eeq
Write $S(u,\l)^{2j-1}$ as a power series of $\l$,
\beq\label{du}
S^{2j-1}(q, \l)=\sum_{i \leq [\frac{2j-1}{2n+1}]+1}S_{2j-1, i}(q)\l^i.
\eeq 
Then coefficients $S_{2j-1, i}(q)$'s are differential polynomials in $q$. 

Assume $j \geq 0$ and $(2j-1) \not\equiv 0 (\rm{mod} (2n+1))$. Then
\beq\label{ak}
q_{t_{2j-1}}=[\p_x+b+q, S_{2j-1, 0}(q)], 
\eeq
is the flow on $C^{\infty}(\R, \cb_n^+)$ constructed from the splitting $(\hat A_{2n}^{(2)})_\pm$ and $J^{2j-1}$. These flows commute. 

Again the group $C^{\infty}(\R, N_n^+ )$ acts on $\p+J+C^{\infty}(\R, \cb_n^+ )$ by gauge transformations with $C^{\infty}(\R, V_n)$ a cross section, and \eqref{ak} is invariant under this action. Hence flow \eqref{ak} induces a quotient flow on $C^{\infty}(\R, V_n)$. We call this quotient flow the \emph{ $(2j-1)$-th $\an 2$-KdV flow}. In fact,  given $u \in C^{\infty}(\R, V_n)$, and $(2j-1) \not\equiv 0\, (\rm{mod} (2n+1))$, there exists a unique differential polynomial $\ti{\eta}_j(u) \in C^{\infty}(\R, \cn_n^+)$ such that
\beq\label{bt}
[\p_x+J+u, S_{2j-1,0}(u)-\ti{\eta}_j(u)] \in C^{\infty}(\R, V_n),
\eeq
where $S_{2j-1,0}(u)$ is defined by \eqref{du}. 

{\it The $(2j-1)$-th $\an2$-KdV flow\/} is
\beq\label{cf}
u_{t_{2j-1}}=[\p_x+b+u, S_{2j-1}(u)-\ti\eta_j(u)],
\eeq
where 
 $S_{2j-1,0}(u)$ is defined as by \eqref{du} and $\ti \eta_j(u)$ is as in \eqref{bt}.

Note that  the $(2j-1)$-th $\an2$-KdV flow can be written as
\beq\label{cfa}
u_{t_{2j-1}}= [\p_x+b+u, P_u(\pi_0(S_{2j-1,0}(u))],
\eeq
where $P_u$ is as defined in Definition \ref{aw} and $\pi_0$ is the projection as in \eqref{bk}.

\bprop\label{dw}\ 

The following statements are equivalent for $u\in C^\infty(\R^2, V_n)$:
\ben 
\item $u$ is a solution of \eqref{cf}.
\item 
\beq\label{ej}
[\p_x+ J+ u, \p_{t_{2j-1}} + (S^{2j-1}(u,\l))_+-\ti\eta_j(u)]=0,
\eeq
for all parameter $\l\in \C$. 
\item 
\beq\label{ek}
[\p_x+b+u, \p_{t_{2j-1}}+ S_{2j-1,0}(u)-\ti \eta_j(u)]=0,
\eeq
\item The linear system 
\beq\label{ep}
\bca g^{-1}g_t= b+u, \\ g^{-1}g_t= S_{2j-1,0}(u)-\ti\eta_j(u),\eca
\eeq is solvable for $g \in C^{\infty}(\R^2, O(n+1, n))$.
\item The linear system 
\beq\label{em}
\bca E^{-1}E_x= J + u, \\ E^{-1}E_t= (S^{2j-1}(u,\l))_+-\ti\eta_j(u),\eca
\eeq is solvable for $E(\cdot, \cdot, \l) \in C^{\infty}(\R^2, O_\C(n+1, n))$ for all parameter $\l\in \C$.
\een
\eprop

We call a solution $E(x,t,\l)$ of \eqref{em} a {\it frame} of the solution $u$ of the $(2j-1)$-th $\an2$-KdV flow if $E(x,t,\l)$ is holomorphic for $\l\in \C$ and  satisfies
$$E(x,t,\l)^{-1}= C_n E^t(x,t,-\l)C_n, \quad \overline{E(x,t,\bar\l)}= E(x,t,\l).$$

Similarly, we have the following.

\bprop\label{cy} If $q: \R^2 \rightarrow \cb_n^+$ is a solution of \eqref{ak}, then $u=\D\ast q$ is a solution of the $(2j-1)$-th $\an2$-KdV flow \eqref{cf}, where $\D(\cdot, t)=D(q(\cdot, t))$ and $D$ is the operator defined in Definition \ref{aba}. Moreover, if $F$ is a frame of the solution $q$ of \eqref{ak}, then $E=F \D^{-1}$ is a frame of the solution $u$ of \eqref{cf}.
\eprop

\bprop\label{ds} Let $u$ be a solution of \eqref{cf}, and $\D:\R^2\to N_n^+$ satisfying $\D_t\D^{-1}=\ti\eta_j(u)$. Then $q(\cdot, t)= \D(\cdot, t)^{-1}\ast u(\cdot, t)$ is a solution of \eqref{ak}.  Moreover, if $E$ is a frame of the solution $u$ of \eqref{cf}, then $F= E\D$ is a frame for the solution $q$ of \eqref{ak}.
\eprop

If $\D_t \D^{-1}=\ti \eta_j(u)$ and $f\in C^\infty(\R, N_n^+)$, then $\D_1(x,t)= \D(x,t) f(x)$ also satisfies $(\D_1)_t\D_1^{-1}=\ti \eta_j(u)$ and $q_1=\D_1\ast u$ is also a solution of \eqref{ak}. 

 Below we write down $S_{2j-1,0}(u)$ and $\ti\eta_j(u)$ for small $n$ and $j$.

\beg\label{eb} {\bf The $\hat A_2^{(2)}$-KdV hierarchy} \

Here $u=q(e_{12}+e_{23})$, and
\begin{align*}
& S_{5, 0}(u)=\bpm -\frac{1}{9}(q_{xxx}-8qq_x) & a & 0  \\ \frac{1}{9}(q_{xx}-4q^2) & 0 & a \\ 0 & \frac{1}{9}(q_{xx}-4q^2) & \frac{1}{9}(q_{xxx}-8qq_x) \epm,  \\
& \ti \eta_{3}(u)= \bpm  0  & \theta & 0 \\ 0 & 0  & \theta \\ 0 & 0 & 0\epm,
\end{align*}
where
\begin{align*}
& a=-\frac{1}{27}q^{(4)}-\frac{7}{27}qq_{xx}+\frac{1}{3}q_x^2+\frac{4}{81}q^3, \\
&\theta=\frac{1}{27}(2q^{(4)}-34qq_{xx}-15q_x^2+\frac{40}{3}q^3).
\end{align*}
Therefore, 
\begin{align*}
&  (S^5(u,\l))_+ - \ti\eta_3(u)\\
&=J^5+\bpm -\frac{q}{3}\l-\frac{1}{9}(q^{(3)}-8qq_x) & -\frac{q_x}{3}\l+\zeta & -\frac{1}{9}(2q_{xx}-5q^2)\l \\
\frac{1}{9}(q_{xx}-4q^2)& \frac{2q}{3}\l & \frac{q_x}{3} \l+\zeta \\
0 & \frac{1}{9}(q_{xx}-4q^2) & -\frac{q}{3}+\frac{1}{9}(q^{(3)}-8qq_x)\l\epm,
\end{align*}
where $\zeta=-\frac{1}{9}(q^{(4)}-9qq_x-8q_x^2+4q^3)$.
It follows that the fifth $\hat A_2^{(2)}$-KdV flow is the KK equation \eqref{kk} given in \cite{K80} and \cite{Ku84}.
\eeg

\beg\label{ga} {\bf The $\hat A_4^{(2)}$-KdV hierarchy} \

We have $u=u_1(e_{23}+e_{35})+u_2(e_{14}+e_{25})$, the third $\hat A_4^{(2)}$-KdV flow is
\beq\label{ac}
\bca
(u_1)_t=-2(u_1)^{(3)}+\frac{7}{2}u_2'+\frac{12}{5}u_1u_1', \\
(u_2)_t=-\frac{5}{4}u_2^{(3)}-\frac{3}{5}(u_1^{(5)}-6u_1u_1^{(3)}+2u_1u_2'-3u_1'u_2).
\eca
\eeq
And 
\begin{align*}
&(S^3(u,\l))_+-\ti\eta_2(u)=\\
&J^3+\bpm -\frac{3}{5}(u_1)_x & \ast & \ast & \ast & \ast \\
-\frac{1}{5}u_1 & -\frac{4}{5}(u_1)_x & \ast & \ast & \ast   \\
0 & \frac{4}{5}u_1 & 0 & \ast & \ast  \\
0 & 0 & \frac{4}{5}u_1  & \frac{4}{5}(u_1)_x & \ast \\
0 & 0 & \ast & -\frac{1}{5}u_1& \frac{3}{5}(u_1)_x \epm.
\end{align*}
\eeg

\ms
Recall that the $j$-th Gelfand-Dickey  (GD$_n$) flow is 
\beq
L_{t_j}=[(L^{\frac{j}{n}})_+, L]
\eeq
for $L=\p^n+\sum_{i=1}^{n-1}v_i\p^{i-1}$, where $(L^{\frac{j}{n}})_+$ is the differential operator component of the pseudo differential operator $L^{\frac{j}{n}}$ (cf. \cite{Dic03}). Next we show that the $\an 2$-KdV hierarchy is a reduction of the GD$_n$ hierarchy (cf. \cite{Dic03}). 

Let $L^\ast$ be the formal adjoint of $L$ under the $L^2$-norm. For example, 
\begin{subequations}
\begin{align}
&\p^*=-\p \\
& f^*=f, \ \ f \in C^{\infty}(\R, \R).
\end{align}
\end{subequations}
It can be checked that
\beq\label{ad}
\cd^s_{2n+1}=\{L=\p^{2n+1}+\sum_{i=1}^{2n}v_i\p^{i-1} \mid L^*=-L\},
\eeq
 is invariant under the $(2j-1)$-th GD$_{2n+1}$ flow, and the induced constraint flow is called the \emph{$(2j-1)$-th constraint KP flow} (cf. \cite{YC92}).  Note that $\cd^s_{2n+1}$ is equal to
\beq
\{L=\p^{2n+1}-\sum_{i=1}^{n}(\p^{n+1-i}u_i\p^{n-i}+\p^{n-i}u_i\p^{n+1-i}) \mid u_i \in C^{\infty}(\R, \R)\}.
\eeq

Consider the matrix eigenvalue problem,
\beq\label{de}
(\p_x+J+u)y=0, 
\eeq
where $J$ is as \eqref{df}, $u=\sum_{i=1}^{n}u_i\b_i$ and $y=(y_1, \cdots, y_{2n+1})^t$. Then the equation for  $y_{2n+1}$ is 
$$(\p^{2n+1}-\sum_{i=1}^{n}(\p^{n+1-i}u_i\p^{n-i}+\p^{n-i}u_i\p^{n+1-i}))y_{2n+1}+\l y_{2n+1}=0.$$

A $\p$-module structure on $\hat A_{2n}^{(1)}$ was introduced in \cite{DS84} to show that $q=\sum_{i=1}^{2n} q_i e_{i, 2n+1}$ is a solution of the $j$-th $\hat A_{2n}^{(1)}$-KdV flow if and only if 
$ \p^{2n+1} -\sum_{i=1}^{2n} q_i \p^{i-1}$
is a solution of the $j$-th GD$_{2n+1}$ flow. It can be checked that $\an2$ is a sub $\p$-module of $\hat A_{2n}^{(1)}$ and a similar proof as in \cite{DS84} gives the following. 

\bprop
$u=\sum_{i=1}^{n}u_i\b_i$ is a solution of the
$(2j-1)$-th $\an2$-flow \eqref{cf} if and only if
$\p^{2n+1}-\sum_{i=1}^{n}(\p^{n+1-i}u_i\p^{n-i}+\p^{n-i}u_i\p^{n+1-i})$ is a solution of the $(2j-1)$-th constraint KP flow.
\eprop

\bs
\section{Hierarchies of isotropic curve flows}\label{hc}

In this section, we  prove result (a) stated in the introduction. 

\bprop
Let $P_{2j-1,0}(u)$ be as in \eqref{dg}, and $S_{2j-1,0}(u)$ as in \eqref{du}. 
Then
 \begin{align}
& \g_t=g P_{2j-1,0}(u)e_1,  \label{cb2} \\
& \g_t= g S_{2j-1,0}(u)e_1.\label{ca2}
\end{align}
are isotropic curve flows on $\cm_{n+1, n}$, where $g(\cdot, t)$ and $u(\cdot, t)$ are the isotropic moving frame and curvature of $\g(\cdot, t)$. We call \eqref{cb} and \eqref{ca}  the {\it $(2j-1)$-th isotropic curve flow on $\cm_{n+1, n}$ of  B-type and  A-type\/} respectively. 
\eprop

\begin{proof} Let $Q_{2j-1}(u)= P_{2j-1,0}(u)-\eta_j(u)$, where $\eta_j(u)$ is as in \eqref{bp}.
 By \eqref{bp}, $[\p_x+b+u, Q_{2j-1}(u)] \in V_n$. It follows from Proposition \ref{bq} that $\xi(\g)=g Q_{2j-1}(u)e_1$ is tangent to $\cm_{n+1, n}$ at $\g$.  Since $\eta_j(u)$ is strictly upper triangular, $Q_{2j-1}(u)e_1= P_{2j-1,0}(u)e_1$. Therefore \eqref{cb} is an isotropic curve flow on  $\cm_{n+1,n}$.  A similar proof shows that \eqref{ca} is an isotropic flow on $\cm_{n+1, n}$.
 \end{proof}
 
Since $\eta_j(u)$ and $\ti\eta_j(u)$ are strictly upper triangular, \eqref{cb2} and \eqref{ca2} can be written as 
 \begin{align}
 &\g_t= g(P_{2j-1,0}(u)-\eta_j(u)) e_1, \label{cb} \\
 & \g_t= g(S_{2j-1,0}(u)-\ti\eta_j(u)) e_1. \label{ca}
 \end{align}
 
\beg {\bf Isotropic curve flows of B-type} \par
\ben
\item[(i)] It follows from Example \ref{el} that the third isotropic curve flow of B-type on $\cm_{n+1, n}$ is \eqref{ar}, i.e., 
$$\g_t=-\frac{3}{n}u_1\g_x+\g_x^{(3)}.$$
In particular,  it is \eqref{bx} when $n=1$, and is \eqref{ap} when $n=2$.
\item[(ii)] By Example \ref{ea}, the fifth isotropic curve flow of B-type on $\cm_{3, 2}$ is
\beq
\g_t=(q_{xxx}-3qq_x)\g+(\frac{3}{2}q^2-q_{xx})\g_x.
\eeq
\een 
\eeg

\beg{\bf Isotropic curve flows of A-Type} \par
\ben
\item[(i)]  The third isotropic curve flow of A-type on $\cm_{n+1, n} (n \geq 3)$ is \eqref{as},
\beq\label{gh}
\g_t=-\frac{3}{2n+1}(u_1)_x\g-\frac{6}{2n+1}u_1\g_x+\g_x^{(3)}.
\eeq
\item[(ii)] By Example \ref{ga},
\beq\label{gg}
\g_t= -\frac{3}{5} (u_1)_x \g - \frac{1}{5} u_1 \g_x + p_4
\eeq is the third isotropic curve flow of A-type on $\cm_{3,2}$. But $p_4= \g_x^{(3)}- u_1 \g_x$ (given in Example \ref{gb}). So \eqref{gg} can be also written as \eqref{gh} with $n=2$. 
\item[(iii)] It follows from the formula of $S_{5,0}(u)$ given in Example \ref{eb} that the  fifth isotropic curve flow of A-type on $\cm_{2, 1}$ is \eqref{al}.
 \een
\eeg

\bthm\label{an} \hfil\par

\ben
\item[(i)] If $\g$ is a solution of the $(2j-1)$-th isotropic curve flow \eqref{cb} of B-type (\eqref{ca} of A-type resp.) on $\cm_{n+1, n}$, then its isotropic curvature $u$ is a solution of the $(2j-1)$-th $\bn1$-KdV flow $\eqref{ai}$ ($\an2$-KdV flow \eqref{cf} resp.).
\item[(ii)] Let $u$ be a solution of \eqref{ai} (\eqref{cf} resp.), $c_0 \in O(n+1, n)$ a constant, and $g$  the solution of  \eqref{eo}
 (\eqref{ep} resp.) with $g(0,0)= c_0$.  Then $\g(x, t):=g(x, t)e_1$ is solution of the $(2j-1)$-th isotropic curve flow of B-type \eqref{cb} (A-type \eqref{ca} resp.) with isotropic curvature $u(\cdot, t)$ and $g(\cdot, t)$ is the isotropic moving frame along $\g(\cdot, t)$. 
\een
\ethm

\begin{proof}
We prove this theorem for \eqref{cb}. The same proof works for \eqref{ca}. 

(i) Let $g(\cdot, t)$ be the isotropic moving frame, and $u(\cdot, t)$ the isotropic curvature along $\g(\cdot, t)$. Then $g^{-1}g_x= b+u$. Set $\xi= g^{-1}g_t$. Hence $b+u$ and $\xi$ satisfy the zero curvature condition, i.e., $[\p_x+b+u, \p_t+\xi]=0$. So we have $u_t= [\p_x+b+u, \xi]$. Since $u_t$ lies in $V_n$, $\xi$ satisfies condition \eqref{ae}.  Note that the first columns of $\xi$ and $Q_{2j-1}(u):=P_{2j-1,0}(u)- \eta_j(u)$ are the same. So by Corollary \ref{gk}, $\xi= Q_{2j-1}(u)$.  By Proposition \ref{dv}, $u$ is a solution of \eqref{ai}.

(ii) Since $g^{-1}g_x= b+u$, $g$ is the isotropic moving frame along $\g(\cdot, t)$. By assumption, $g^{-1}g_t= Q_{2j-1}(u)$. Hence $g_t= gQ_{2j-1}(u)$. Note that the first column of $g$ is $\g$. So we have $\g_t= gQ_{2j-1}(u)e_1$. 
\end{proof}

\bcor\label{et} Let $\Psi:\cm_{n+1,n}\to C^\infty(\R, V_n)$ be the isotropic curvature map. Then $\Psi$ maps the space of solutions of \eqref{cb}  (\eqref{ca} respectively) modulo $O(n+1, n)$  bijectively onto the space of solutions of \eqref{ai} (\eqref{cf} respectively).
\ecor

\ms

The following Theorem is a consequence of Theorem \ref{an}.

\bthm {\bf [Cauchy problem on the line]} \ 

\ni Let $\g_0 \in \cm_{n+1, n}$ with rapidly decaying isotropic curvature $u_0$, and $g_0 \in O(n+1, n)$ the isotropic moving frame alogn $\g_0$. Let $u(x, t)$ be the solution of the $(2j-1)$-th $\bn 1$-KdV flow \eqref{ai} ($\an 2$-KdV flow \eqref{cf} resp.) with initial date $u(x, 0)=u_0(x)$, and $g$ the solution of \eqref{eo} (\eqref{ep} resp.) with $g(0,0)=g_0(0)$. Then $\g(x, t)=g(x,t)e_1$ is a solution of \eqref{cb} (\eqref{ca} respectively) with $\g(x, 0)=\g_0(x)$ and $\g(\cdot, t)$ has rapidly decaying isotropic curvatures.
\ethm

If the solution of the periodic Cauchy problem for  \eqref{ai} is solved, then we can use a similar proof as for the $\hat A_n^{(1)}$-KdV flows (cf. \cite{TWb}) to solve the periodic Cauchy problem for \eqref{eo}, i.e., we have the following.

\bthm {\bf [Cauchy problem with periodic initial data]} \ 

\ni Suppose $\g_0\in \cm_{n+1, n}$ is periodic, $g_0$ and $u_0$ are the isotropic moving frame and curvature along $\g_0$. Let $u(x, t)$ be the solution of the $(2j-1)$-th $\bn 1$-KdV flow \eqref{ai} ($\an 2$-KdV flow \eqref{cf} resp.) periodic in $x$ such that $u(x, 0)=u_0(x)$, and $g(x,t)$ the solution of \eqref{eo} (\eqref{ep} resp.) with $g(0,0)=g_0(0)$. Then $\g(x, t)=g(x,t)e_1$ is a solution of \eqref{cb} (\eqref{ca} resp.) with $\g(x, 0)=\g_0(x)$. Moreover, $\g(\cdot , t)$ is periodic in $x$ with isotropic curvature $u(\cdot, t)$.
\ethm

\bigskip
\section{Bi-Hamiltonian structure for isotropic curve flows}\label{he}

In this section, we first explain how to compute the bi-Hamiltonian structure and conservation laws for the $\bn1$-KdV and $\an2$-KdV hierarchies. Then we  pull back these bi-Hamiltonian structures and conservation laws by the isotropic curvature map $\Psi$ to get bi-Hamiltonian structures and conservation laws for isotropic curve flows. 

The gradient $\K F(u)\in C^\infty(S^1, V_n^t)$ for a functional $F$ on $C^\infty(S^1, V_n)$ is defined by 
$$\rd F_u(v)= \li \K F(u), v\ri =\oint_{S^1} \tr(v\K F(u))\rd x$$
for all $v\in C^\infty(S^1, V_n)$. 

If $\{\, ,\}$ is a Poisson structure on $C^\infty(S^1, V_n)$, then the Hamiltonian vector field $X_F$ for $F$ with respect to $\{\, ,\}$ is defined by 
\beq\label{eu}
\{ F, H\}(u) =-\li X_F(u),  \K H(u)\ri
\eeq
for all functionals $H$. 

The bi-Hamiltonian structure on $C^\infty(S^1, V_n)$ for the $\bn 1$-KdV hierarchy given in \cite{DS84} can be written as follows:
\begin{align*}
& \{F_1, F_2\}_1(u)=\li [\b, P_u(\K F_1(u))], P_u(\K F_2(u))\ri, \\
& \{F_1, F_2\}_2 (u)=\li [\p_x+b+u, P_u(\K F_1(u))], P_u(\K F_2(u))\ri,
\end{align*}
where $\b$ is defined by \eqref{db} and $P_u: C^{\infty}(S^1, V_n^t) \rightarrow C^{\infty}(S^1, o(n+1, n))$ is the linear operator defined in Definition \ref{aw}.  Use \eqref{eu} to see that the Hamiltonian vector field $Y_F$ of a functional $F$ with respect to $\{\, ,\}_2$ is  
\beq\label{gi}
Y_F(u)= [\p_x+b+u, P_u(\K F(u))].
\eeq
We have explained how to compute $P_u(\xi)$ in section \ref{hf}. So we can compute the Hamiltonian vector field $X_F$ of $F$ with respect $\{\, ,\}_1$.

\beg {\bf Bi-Hamiltonian structure for the $\hat B_1^{(1)}$-KdV hierarchy}\ \\
Write $\ti{\xi}=\K F_1(u)=\xi(e_{21}+e_{32})$, $\ti{\eta}=\eta(e_{21}+e_{32})$, $C=P_u(\ti{\xi})=(C_{ij})$ and $D=P_u(\eta)=(D_{ij})$. We use Example \ref{dk} to get $C$ and $D$ in terms of $\xi$ and $\eta$ respectively. So we have 
\begin{align*}
&\{F_1, F_2\}_1(u)=\li [C, \b], D \ri=-2\oint \xi_x\eta dx, \\
&\{F_1, F_2\}_2(u)=\li [\p_x+b+u, C], D\ri=-2\oint (\xi_{xxx}-2u_1\xi_x-(u_1)_x\xi)\eta dx.
\end{align*}
This is the standard bi-Hamiltonian structure for the KdV-hierarchy \cite{Dic03}.
\eeg

\beg {\bf Bi-Hamiltonian structure for the $\hat B_2^{(1)}$-KdV hierarchy}\ \\
Write ${\xi}=\xi_1(e_{32}+e_{43})+\xi_2(e_{41}+e_{52})$, ${\eta}=\eta_1(e_{32}+e_{43})+\eta_2(e_{41}+e_{52})$, $C=P_u(\xi)=(C_{ij})$ and $D=P_u(\eta)=(D_{ij})$. We use Example \ref{dl} to get $C$ and $D$ in terms of $ \xi$ and $ \eta$ respectively, and obtain
\begin{align*}
&\{F_1, F_2\}_1(u)=\li [C, \beta_2], D \ri \\
&  \quad =(C_{11}+C_{22})\eta_2+C_{31}\eta_1-\xi_1D_{31}-\xi_2(D_{11}+D_{22}) \\
& \quad=-2\oint((\xi_2)_x^{(3)}+2(\xi_1)_x-(u_1\xi_2)_x-u_1(\xi_2)_x)\eta_2+2(\xi_2)_x\eta_1 dx, \\
&\{F_1, F_2\}_2(u)=\li [\p_x+b+u, C], D \ri=2\oint p_1\eta_2+p_2\eta_1 dx,
\end{align*}
where $p_1=[\p_x+b+u, C]_{14}$ is a $7$-th ordered differential polynomial in $\xi_1, \xi_2$, and $p_2=[\p_x+b+u, C]_{23}$ is a $5$-th ordered differential polynomial in $\xi_1. \xi_2$.
\eeg

The following theorem can be proved the same way as for the $\hat A_n^{(1)}$-KdV hierarchy (cf. \cite{TWb}).

\bthm\label{co} Let $u \in \C^{\infty}(\R, V_n)$, $\b$ as in \eqref{db}, and $P(u,\l)$ defined by \eqref{ch}.  Then we have
\beq\label{er}
\li \frac{\p }{\p\l} (\l^{-1}P^{2j-1}(u, \l)), \d u\ri= \d \li \l^{-1} P^{2j-1}(u, \l), \b \ri=  \li \d P^{2j-1}(u,\l),  \b \l^{-1} \ri.
\eeq
\ethm

\bthm\label{cp} Let $u, \b, P(u,\l)$ be as in Theorem \ref{co}. Write $P^{2j-1}(u,\l)= \sum_i P_{2j-1, i}(u)\l^i$ as a power series in $\l$. Set 
\beq\label{cq}
F_{2j-1}(u)=-\oint \tr(P_{2j-1, -1}(u)\b) \rd x.
\eeq
Then $\K F_{2j-1}(u)=\pi_0(P_{2j-1,0}(u))$, where $\pi_0$ is the projection onto $V_n^t$ defined by \eqref{bk}.
 Moreover, the Hamiltonian equation for $F_{2j-1}$ with respect to $\{\, ,\}_2$ ($\{\, ,\}_1$ resp.) is the $(2j-1)$-th ($(2(j-n)-1)$-th resp.) $\bn1$-KdV flow.
 \ethm

\begin{proof} Compare the coefficient of $\l^{-2}$ of \eqref{er} to obtain the formula for $\K F_{2j-1}$. 

By Corollary \ref{cr}, we have 
$P_u(\pi_0(P_{2j-1, 0}(u))=  P_{2j-1,0}(u)-\eta_j(u)$. So 
\beq\label{dm}
P_u(\K F_{2j-1}(u))= P_{2j-1,0}(u,\l)-\eta_j(u).
\eeq 
It follows from \eqref{gi} that the Hamiltonian flow for $F_{2j-1}$ with respect to $\{\, \}_2$ is the $(2j-1)$-th $\bn1$-KdV flow. 

Let $X_{2n+2j-1}$ denote the Hamiltonian vector field of $F_{2n+2j-1}$ with respect to $\{\, ,\}_1$. Compute directly to get 
\begin{align*}
& \{F, F_{2n+2j-1}\}_1(u) = -\li X_{2n+2j-1}(u), \K F(u)\ri\\
&=-\li [\b, P_u(\K F_{2n+2j-1}(u))], P_u(\K F(u))\ri, \, {\rm by\,\, \ref{dm},}\\
&= -\li [\b, P_{2n+2j-1,0}(u)-\eta_{n+j}(u)], P_u(\K F(u))\ri, \quad {\rm since\,} [\b, \eta_j(u)]=0,\\
&=-\li [\b, P_{2n+2j-1,0}(u)], P_u(\K F(u))\ri.
\end{align*}
Compare coefficient of $\l$ of the equation \eqref{ch} to get
$$[\p_x+b+u, P_{2n+2j-1,1}(u)]= [P_{2n+2j-1,0}(u), \b].$$
It follows from \eqref{dh} that we have
\beq\label{dn}
P^{2n+2j-1}(u,\l)= \l P^{2j-1}(u,\l).
\eeq 
Compare coefficient $\l$ of \eqref{dn} to get $P_{2n+2j-1, 1}(u)= P_{2j-1, 0}(u)$. So we have
$$\{F, F_{2n+2j-1}\}_1(u)= \li [\p_x+b+u, P_{2j-1,0}(u)], P_u(\K F(u))\ri.$$
Since $\eta_j(u)\in \cn_n^+$ and $[\p_x+b+u, P_u(\K F(u))]\in V_n$, 
$$\li [ \p_x+b+u, \eta_j(u)], P_u(\K F(u))\ri = -\li \eta_j(u), [\p_x+b+u, P_u(\K F(u))]=0.$$
This implies that 
$$\{F, F_{2n+2j-1}\}_1(u)= \li [\p_x+ b+u, P_u(\K F_{2j-1}(u))], P_u(\K F(u))]\ri.$$
By the definition of $P_u$, we have $[\p_x+b+u, P_u(\K F_{2j-1}(u))]\in V_n$ and $\pi_0(P_u(\K F(u))= \K F(u)$. Hence 
$$ \{F, F_{2n+2j-1}\}_1(u)= \li [\p_x+b+u, P_u(\K F_{2j-1}(u))], \K F(u)\ri.$$
This proves that the Hamiltonian flow for $F_{2(j+n)-1}$ is the $(2j-1)$-th $\bn1$-KdV flow.  
\end{proof}

\beg {\bf Conservation laws for the $\bn1$-KdV hierarchy}\

Let $h_{2j-1}(u)=\tr(P_{2j-1, -1}(u)\b)$ denote the density of the conservation law $F_{2j-1}$.  
\ben
\item For $n=1$,  we have  $u=q(e_{12}+e_{23})$ and
$$h_1=2q, \ \ h_3=q^2, \ \ h_5=\frac{1}{2}(q^3-qq_{xx}).$$
\item For $n=2$, we have $u=u_1(e_{23}+e_{34})+u_2(e_{14}+e_{25})$ and
$$h_1=\frac{1}{2}u_1, \ \ h_3=\frac{1}{8}u_1^2+\frac{1}{2}u_2.$$
\item For general $n$, we have
$$h_1=\frac{1}{n}u_1, \quad h_3=\frac{2n-3}{2n^2}u_1^2+\frac{1}{n}u_2. $$
\een
\eeg

We use the same proofs for the $\hat A_n^{(1)}$-KdV hierarchy to prove the following results for the $\an2$-KdV. 

\bthm\label{gj} Let $u \in C^{\infty}(\R, V_n)$, and $S(u, \l)$ defined by \eqref{du}. Then 
\begin{equation*}
\li \frac{\p }{\p \l}(\l^{-1}S^{2j-1}(u, \l)), \d u \ri=\l^{-1}\li \d S^{2j-1}(u, \l), e_{1, 2n+1} \ri .
\end{equation*}
\ethm

\bthm Let $u, S(u,\l)$ be as in Theorem \ref{gj}, and $S^{2j-1}(u, \l)=\sum_{i}S_{2j-1, i}(u)\l^i$. Set $h_{2j-1}(u)=  tr(S_{2j-1, -1}(u)e_{1, 2n+1})$, and
\beq
G_{2j-1}(u)=-\oint h_{2j-1}(u)\rd x.
\eeq
Then $\nabla G_{2j-1}(u)=\pi_0(S_{2j-1, 0}(u))$, where $\pi_0$ is the projection defined by \eqref{bk}. Moreover, the Hamiltonian flow for $G_{2j-1}$ with respect to $\{\, ,\}_2$ ($\{\, ,\}_1$ resp.) is the $(2j-1)$-th ($(2(j-n)-1)$-th resp.) $\an2$-KdV flow.
\ethm

\beg {\bf Conservation laws for the $\an2$-KdV hierarchy}\ 
\ben
\item For the $\hat A_{2}^{(2)}$-KdV hierarchy, $u=q(e_{12}+e_{23})$, we have 
\begin{equation*}
h_1(u)=\frac{2}{3}q, \ \ h_3(u)=\frac{5}{27}(q_x^2+\frac{8}{3}q^3).
\end{equation*}
\item For general $n$, the first two conservation densities are
\begin{align*}
& h_1(u)= \frac{2}{2n+1}u_1, \quad h_3(u)= \frac{2}{2n+1}u_2+\frac{4(n-1)}{(2n+1)^2}u_1^2 \end{align*}
\een
\eeg

Let $\cm_{n+1, n}(S^1)$ denote the space of $\g\in \cm_{n+1, n}$ that is periodic with period $2\pi$. 
Then we have the following.
\ben
\item[(i)] The isotropic curvature map $\Psi:\cm_{n+1, n}(S^1)\to C^\infty(S^1, V_n)$ induces an injective map from the orbit space $\frac{\cm_{n+1, n}(S^1)}{O(n+1, n)}$ to $C^\infty(S^1, V_n)$.
\item[(ii)] The isotropic curve flows are invariant under $O(n+1, n)$.
\item[(iii)] Suppose $\g_t= g\xi(u)$ is an isotropic flow on $\cm_{n+1,n}$. By Proposition \ref{ao}, its curvature evolves as $u_t= [\p_x+b+u, g^{-1}g_t]$, where $g(\cdot, t)$ is the isotropic moving frame. Moreover, it follows from Corollaries \ref{eh} and \ref{gk} that $g^{-1}g_t$ can be computed from $\xi(u)$.  
\een

We pull back the bi-Hamiltonian structure and conservation laws for the $\bn1$-KdV and $\an2$-KdV hierarchies on $C^\infty(S^1,V_n)$ to $\cm_{n+1, n}(S^1)$ by the isotropic curvature map $\Psi$ for the isotropic curve flows. In particular, we have the following:
 \ben
 \item[(a)] A functional on $\frac{\cm_{n+1, n}(S^1)}{O(n+1, n)}$ can be viewed as an $O(n+1, n)$-invariant functional on $\cm_{n+1,n}(S^1)$. So it is of the form 
 $\hat F(\g)= F(\Psi(\g))$ for some functional $F$ on $C^\infty(S^1, V_n)$.
 \item[(b)] Given functionals $F, H$ on $C^\infty(S^1, V_n)$, let $\hat F= F\circ \Psi$ and $\hat H= H\circ \Psi$. Then the induced Poisson structure on $O(n+1, n)$-invariant functionals on $\cm_{n+1, n}(S^1)$ is 
$$ \{\hat F, \hat H\}^\wedge_i(\g)= \{ F, H\}_i(\Psi(u)), \quad i=1, 2.$$
\een

The following is a consequence of Proposition \ref{ao}, Corollary \ref{es}, and Corollary \ref{eh}. 

\bthm Let $\Psi$ be the isotropic curvature map, and  $X_F$ the Hamiltonian vector field of $F:C^\infty(S^1, V_n)\to \R$ with respect to $\{\, ,\}_1$. Then 
\ben
\item[(i)] there exists $\xi_F(u)\in C^\infty(S^1, o(n+1, n))$ satisfying 
$$[\p_x+b+u, \xi_F(u)]= X_F(u),$$ 
\item[(ii)] the Hamiltonian equation for $\ti F= F\circ \Psi$ on $\cm_{n+1, n}$ with respect to $\{\, ,\}_1^\wedge$ ($\{\, ,\}^\wedge_2$ resp.) is $\g_t= g\xi_F(u) e_1$ ($\g_t= gP_u(\K F(u)) e_1$ resp.),
\een
where $g(\cdot, t)$ and $u(\cdot, t)$ are  the isotropic curvature frame and isotropic curvature along $\g(\cdot, t)$. 
\ethm

In particular, the Hamiltonian flow for $\hat F_{2j-1}$ with respect to $\{\, ,\}^\wedge_2$ ($\{\, ,\}_1^\wedge$ resp.) on $\cm_{n+1, n}$ is the $(2j-1)$-th ($(2(j-n)-1)$-th resp.) isotropic curve flow of B-type. Similar statements hold for isotropic curve flow of A-type.

\bs

\section{B\"acklund transformations for the $\bn1$-KdV flows}\label{hd}

In this section, we first construct B\"acklund transformations (BTs) and a Permutability formula for flow \eqref{ah} on $C^\infty(\R, \cb_n^+)$. Then we use the gauge equivalence to construct BTs for the corresponding quotient flow \eqref{ai}, i.e., the $(2j-1)$-th $\bn1$-KdV flow. Since we also obtain the formula of the frame of the new solution constructed from BTs for \eqref{ah}, we can construct BTs for isotropic curve flows of B-type. If we apply BTs to the trivial solution of the isotropic curve flow (i.e., the solution whose isotropic curvatures are zero) repeatedly, then we can obtain infinitely many families of explicit soliton solutions.

Let $\hat\B_n^{(1)}$ denote the group of smooth maps $f:S^1\to SL(2n+1, \C)$ satisfying 
\beq\label{bu}
\bca
\overline{f(\bar \l)}=f(\l), \\
f(\l)^tC_nf(\l)=C_n,
\eca
\eeq
and  $(\hat \B_n^{(1)})_+$  the subgroup of $f\in \hat \B_n^{(1)}$ that is the restriction of a holomorphic map on $\C$ to $S^1$, and $(\hat \B_n^{(1)})_-$ the subgroup of $f\in \hat \B_n^{(1)}$ that is the boundary value of a holomorphic map $\ti f$ on $\e^{-1}\leq |\l |\leq \infty$ for some small $\e>0$ and $\ti f(\infty)=\I $. Then the Lie algebras of $\hat \B_n^{(1)}$ and $(\hat \B_n^{(1)})_\pm$ are $\bn1$ and $(\bn1)_\pm$ respectively. If a soliton hierarchy is constructed from a splitting $\cl_\pm$ of a loop algebra $\cl$, then we can use the loop group factorization to constructing BTs (cf. \cite{TU00}) as follows:
\ben
\item Find {\it simple elements} ( i.e., rational maps) $f\in (\bbn1)_-$ that have minimum number of poles.
\item Given $f\in (\bbn1)_-$ and a frame $F(x,t,\cdot)\in (\bbn1)_+$ of a solution of \eqref{ah}, if we can factor $fF(x, t, \cdot)=\ti F(x, t, \cdot)\ti f(x, t, \cdot)$ with $\ti F(x,t,\cdot)\in (\bbn1)_+$ and $\ti f(x, t, \cdot) \in (\bbn1)_-$, then it was proved in \cite{TU00} that $\ti F$ is a frame of a new solution of \eqref{ah}.
\een

We need the following Lemmas to construct rational elements in $(\bbn1)_-$.

\blem\label{cs} Let $\R^{n+1, n}= V_1\oplus V_2$, and $V_i^\perp=\{v\in \R^{n+1, n}\n  \li v, V_i\ri =0\}$.   Let $\pi$ be the projection of $\R^{n+1, n}$ onto $V_1$ along $V_2$, and $\pi^\sharp$  the projection of $\R^{n+1, n}$ onto $V_2^{\perp}$ along $V_1^\perp$. Then $\pi^{\sharp}=C_n\pi^tC_n$.
\elem
\begin{proof} Note that $\pi^{\sharp}$ is the conjugate of $\pi$ with respect to $\li \ , \ \ri$. First, if $\li \pi^{\sharp} X, Y \ri = 0$ for all $X \in \R^{n+1, n}$, then $\li X, \pi Y \ri=0$. Hence $\pi Y = {\bf 0}$. Therefore, $\Im(\pi^{\sharp})^\perp \subset \Ker(\pi)=V_2$. Since $\pi^{\sharp}=C_n\pi^tC_n$ and $\dim(\Im(\pi^{\sharp}))=\dim(\Im(\pi))$, we have  $\Im(\pi^{\sharp})^\perp=V_2$. Hence $\Im(\pi^{\sharp})=V_2^{\perp}$.

If $Y \in \Ker(\pi^{\sharp})$, then $\li \pi X, Y \ri=\li X, \pi^{\sharp}Y \ri=0$ for all $X \in \R^{n+1, n}$. Hence $Y \in V_1^{\perp}$. Therefore, $\Ker(\pi^{\sharp})=V_1^{\perp}$.
\end{proof}

\blem \label{bc}  Let $V_i, \pi, \pi^\sharp$ be as in Lemma \ref{cs}.  Then $\pi \pi^{\sharp}=\pi^{\sharp}\pi={\bf 0}$ if and only if
\beq\label{cu}
V_1 \subset V_1^{\perp}, \quad V_2^{\perp} \subset V_2.
\eeq
\elem

\begin{proof} Since $V_2^{\perp} \subset V_2$, $\pi \pi^{\sharp}=\bf{0}$. It follows from $V_1 \subset V_1^{\perp}$ that we have $\pi^{\sharp}\pi=\bf{0}$.
\end{proof}

 Let $\a_1 \neq \a_2 \in \R$, and $\pi$ a projection of $\R^{n+1, n}$. Set
\beq\label{ct}
h_{\a_1, \a_2, \pi}=I+\frac{\a_1-\a_2}{\l-\a_1}(I-\pi) = \frac{\l-\a_2}{\l-\a_1}\I -\frac{\a_1-\a_2}{\l-\a_1} \pi.
\eeq

\bprop\label{bca}
If $\pi\pi^\sharp= \pi^\sharp\pi=0$, then
\beq\label{bd}
g_{\a_1, \a_2, \pi}=h_{\a_2, \a_1, \pi^{\sharp}}h_{\a_1, \a_2, \pi}.
\eeq
is in $(\bbn1)_-$.
\eprop

\begin{proof} By Lemma \ref{bc}, we have $\pi\pi^\sharp=\pi^\sharp \pi=0$. So we obtain
\begin{align*}
&g_{\a_1, \a_2, \pi}=I+\frac{\a_2-\a_1}{\l -\a_2}\pi+\frac{\a_1-\a_2}{\l-\a_1}\pi^{\sharp},\\
&g_{\a_1, \a_2, \pi}^{-1}=I+\frac{\a_2-\a_1}{\l -\a_2}\pi^{\sharp}+\frac{\a_1-\a_2}{\l-\a_1}\pi=g_{\a_1, \a_2, \pi^{\sharp}}.
\end{align*}
A direct computation implies that  $C_ng_{\a_1, \a_2, \pi}^{-1}= g_{\a_1, \a_2, \pi}^tC_n$.
\end{proof}

\bthm \label{bea}
Let $F(x, t, \cdot) \in (\bbn1)_+$ be a frame of a solution $q$ of \eqref{ah}. Assume that $\R^{n+1,n}=V_1\oplus V_2$ and $V_1, V_2$ satisfy \eqref{cu}.  Let  $\pi$ be the projection of $\R^{n+1,n}$ onto $V_1$ along $V_2$, and $\pi^{\sharp}=C_n \pi^t C_n$. Let $\a_1 \neq \a_2 \in \R$, $g_{\a_1, \a_2, \pi}$  as in \eqref{bd}, $\tilde V_i(x,t)=F(x, t, \a_i)^{-1}V_i, i=1,2$,  $\ti\pi(x,t)$ the projection of $\R^{n+1, n}$ onto $\ti V_1(x,t)$ along $\ti V_2(x,t)$, and  
$$\ti{F}(x, t, \l)=g_{\a_1, \a_2, \pi}(\l)F(x,t,\l)g^{-1}_{\a_1, \a_2, \ti\pi(x,t)}(\l).$$
Then 
$$\ti{q}=q+(\a_1-\a_2)[\b, \ti\pi-\ti\pi^{\sharp}]$$ is a new solution of \eqref{ah} and $\ti F(x,t,\cdot)\in (\bbn1)_+$ 
  is a frame for $\ti q$ .
 (We will use $g_{\a_1, \a_2, \pi}\sharp q$ to denote the new solution $\ti q$). 
 \ethm

\begin{proof} Since $F(x, t, \cdot) \in (\bbn 1)_+$ and $\a_i\in \R$, $F(x, t, \a_i)\in O(n+1, n)$. So 
\begin{align*}
\li \ti{V}_1, \ti{V}_1 \ri=\li F(x, t, \a_1)^{-1}V_1,  F(x, t, \a_1)^{-1}V_1 \ri= \li V_1, V_1\ri=0.
\end{align*}
Therefore, $\ti{V}_1 \subset \ti{V}_1^{\perp}$. 

  If $v \in \ti{V}_2^{\perp}$, then 
\begin{align*}
\li F(x,t, \a_2) v, Y\ri=\li v, F(x, t, \a_2)^{-1}Y \ri=0, \quad  \forall  \ Y \in V_2.
\end{align*}
Hence $F(x,t, \a_2) v \in V_2^{\perp}$. This implies that
$$\ti{V}_2^{\perp} \subset F(x,t, \a_2)^{-1}(V_2^{\perp}) \subset F(x,t, \a_2)^{-1}(V_2)=\ti{V}_2.$$
By Lemma \ref{bc} and Proposition \ref{bca}, $g_{\a_1, \a_2, \ti{\pi}(x,t)}\in (\bbn1)_-$.

Next we claim that $\ti{F}(x, t, \l)$ is holomorphic for $\l \in \C$. It follows from the formula of $\ti{F}$ that $\ti F(x, t, \l)$ is holomorphic for $\l \neq \a_1, \a_2$. The residue of $\ti{F}(x, t, \l)$ at $\a_1$ is: 
$$(\a_1-\a_2)(\pi^{\sharp}F(x, t, \a_1)(I-\ti{\pi}^{\sharp})+(I-\pi)F(x, t, \a_1)\ti{\pi}).$$
Since $\ti{\pi}^{\sharp}$ is the projection of $\R^{n+1, n}$ onto $\ti{V}_2^{\perp}$ along $\ti{V}_1^{\perp}$, we have
$$\Im(F(x, t, \a_1)(I-\ti{\pi}^{\sharp})) \subset V_1^{\perp}, \quad \Im(F(x, t, \a_1)\ti{\pi}) \subset V_1.$$
Hence $\pi^{\sharp}F(x, t, \a_1)(I-\ti{\pi}^{\sharp})+(I-\pi)F(x, t, \a_1)\ti{\pi}=0$. So $\ti F(x, t, \l)$ is holomorphic at $\l=\a_1$.

A similar computation implies that the residue of $\ti{F}(x,t, \l)$ at $\a_2$, 
$$(\a_2-\a_1)(\pi F(x, t, \a_2)(I-\ti{\pi})+(I-\pi^{\sharp})F(x, t, \a_2)\ti{\pi}^{\sharp})=0.$$
So we have proved that $g_{\a_1, \a_2, \pi}F(x, t, \cdot)=\ti F(x, t, \cdot)g_{\a_1, \a_2, \ti \pi(x, t)}$
with $\ti F(x, t, \cdot) $ in $ (\bbn1)_+$ and $g_{\a_1, \a_2, \ti \pi(x, t)}$ in $(\bbn1)_-$. It follows from \cite{TU00} that $\ti F$ is a frame of a new solution $\ti q$ of \eqref{ah}. 

Note that $\ti F^{-1}\ti F_x= b+\ti q+\l \b$, $F^{-1}F_x= b+q+\l\b$, and 
$$\ti F^{-1}\ti F_x= g_{\a_1, \a_2, \ti\pi} F^{-1}F_x  g_{\a_1, \a_2, \ti\pi}^{-1} - ( g_{\a_1, \a_2, \ti\pi})_x g_{\a_1, \a_2, \ti\pi}^{-1}.$$
So we have
$$(b+\ti q+\b\l)  g_{\a_1, \a_2, \ti\pi}=  g_{\a_1, \a_2, \ti\pi} (b+ q+\b\l) - ( g_{\a_1, \a_2, \ti\pi})_x.$$
Equate the constant term of the power series expansion of the above equation to get the formula for $\ti q$.
\end{proof}

As a consequence of Proposition \ref{cga}, Proposition \ref{do}, and Theorem \ref{bea}, we obtain BTs for the $(2j-1)$-th $\bn1$-KdV flow.

\bthm\label{be} Let $E(x,t,\cdot)\in (\bbn1)_+$ be a frame of a solution $u$ of \eqref{ai}, $u_t= [\p_x+b+u, P_{2j-1,0}(u)-\eta_j(u)]$, and $\D(x,t)\in N_n^+$ satisfying $\D_t\D^{-1}= \eta_j(u)$. Let $q=\D^{-1}\ast u$,  $\a_1, \a_2, \pi, \ti \pi$  as in Theorem \ref{bea}, and $\ti \D=D(\ti q)$, where $\ast$ is the action defined by \eqref{dr} and $D$ is the operator given in Definition \ref{aba}.  Then 
\ben
\item[(i)] $q$ is a solution of \eqref{ah},
\item[(ii)]  let $\ti q= g_{\a_1, \a_2, \ti\pi}\sharp q$ be as in Theorem \ref{bea}, and $$\ti E= g_{\a_1, \a_2, \pi}E\D g_{\a_1,\a_2, \ti\pi}^{-1}\ti \D^{-1}.$$  Then $\ti u= \ti \D\ast(g_{\a_1, \a_2, \pi^\sharp} (\D^{-1} \ast u))$ is a solution of \eqref{ai} and $\ti E$ is a frame of $\ti u$. 
\een
 \ethm
 
\bthm\label{ci} Let $\g: \R^2 \rightarrow \cm_{n+1, n}$ be a solution of the $(2j-1)$-th isotropic curve flow  \eqref{cb} of  B-type, $g(\cdot, t)$ the isotropic moving frame along $\g(\cdot, t)$, and $u(\cdot, t)=(g^{-1}g_x-b)$ the solution of the $(2j-1)$-th $\bn 1$-KdV flow \eqref{ai} as in Theorem \ref{an}. Let $\a_1, \a_2, \pi, \pi^{\sharp}, \ti \pi, \ti \pi^{\sharp}, \D, \ti u$ and $\ti E$ be as in Theorem \ref{be}. Then $\ti \g(x,t)= \ti E(x,t,0)e_1$ is a solution of \eqref{cb} with isotropic curvature $\ti u$, where $e_1=(1, 0, \cdots, 0)^t$.
\ethm

\begin{proof}
Let $E(x,t,\l)$ be the frame of the solution $u$ of \eqref{ai} satisfying $E(0,0,\l)=g(0,0)$. Note that $E(x,t,0)$ and $g(x,t)$ satisfy the same linear system,
$$\bca g^{-1}g_x= b+ u, \\ g^{-1}g_t= P_{2j-1,0}(u)-\eta_j(u),\eca$$
and have the same initial data. So $E(x,t,0)= g(x,t)$ for all $x, t$.  By Theorem \ref{be}, $\ti E$ is a frame of a new solution $\ti u$. It follows from Theorem \ref{an} (ii) that $\ti \g$ is a solution of \eqref{cb}. 
\end{proof}

Next we write down the formula of BT for $\g$ in terms of $\g$.

\bcor\label{dx} Let $\g, u, g, \a_1, \a_2, \pi, \D, \ti\pi$ be as in Theorem \ref{ci}.
 If  $\a_1\a_2\not=0$, then
$$
\hat \g(x, t)=g(x, t)\D(x,t)\left(I+\frac{\a_1-\a_2}{\a_2}\ti \pi^{\sharp}(x, t)+\frac{\a_2-\a_1}{\a_1}\ti \pi(x,t)\right)e_1$$
is a solution of \eqref{cb}.
\ecor 

\begin{proof}  By Theorem \ref{ci}, $\ti \g(x,t)=\ti E(x,t,0) e_1$ is a solution of the \eqref{cb}.  Note that $g_{\a_1, \a_2, \pi}(\l)$ is holomorphic at $\l=0$ and $g_{\a_1, \a_2, \pi}(0)\in O(n+1, n)$. So  
$$\hat \g= E(x,t,0)\D(x,t) g^{-1}_{\a_1, \a_2, \ti \pi(x,t)}(0) \ti \D(x,t)^{-1}e_1$$
is also a solution of \eqref{cb}. 
  But $\ti \D \in N_n^+$ implies that $\ti \D(x,t)e_1= e_1$. This proves the Corollary. 
  \end{proof}
  
\brem 
 If $\a_1\not=0$ and $\a_2=0$, the formula for the new solution $\ti \g$ obtained in Theorem \ref{ci} is more complicated. 
We use the same notation as in Theorem \ref{ci} and set
$$\pi'=\I-\pi, \quad (\pi^{\sharp})'=\I-\pi^{\sharp}.$$ Write $E(x, t, \l)$ as a power series in $\l$,
$$E(x, t, \l)=E_0(x, t)+E_{1}(x, t)\l+E_2(x, t)\l^2 + \cdots.$$ Then the new solution $\ti \g$ obtained in Theorem \ref{ci} is\begin{align*}
\ti \g(x, t)& = ( (\pi^{\sharp})'(E_0(x, t)\D \ti \pi'-\a_1E_1(x, t)\D \ti \pi^{\sharp}) e_1 \\
& \quad +(-\a_1\pi E_1(x, t)\D \ti \pi'+\a_1^2 \pi E_2(x, t)\D \ti \pi^{\sharp}) e_1.
\end{align*}
\erem

\beg {\bf $1$-soliton of the third isotropic curve flow of $B$-type on $\cm_{2,1}$}\ 

Note that $\g=(1, x, \frac{x^2}{2})^t$ is a solution of the third isotropic curve flow \eqref{bx}  with isotropic curvature $u=0$ and $E(x,t,\l)= \exp(J_Bx + J_B^3t)$ is a frame of $u=0$.   Next we apply BTs  to write down explicit solutions for the  isotropic curve flow \eqref{bx}, 
$$\g_t= u_x \g- u\g_x.$$  Note that $\eta_3(0)=0$, so we can choose $\D$ to be the identity.

Set $\l=z^2$. Then $E(x, t, z^2)$ is equal to 
$$
\bpm  \frac{1}{2}(\cosh(zx+z^3t)+1) & \frac{z}{2}\sinh(zx+z^3t) & \frac{z^2}{4}(\cosh(zx+z^3t)-1) \\
\frac{1}{z} \sinh(zx+z^3t) & \cosh(zx+z^3t) &  \frac{z}{2}\sinh(zx+z^3t) \\
\frac{1}{z^2}(\cosh(zx+z^3t)-1) & \frac{1}{z} \sinh(zx+z^3t) & \frac{1}{2}(\cosh(zx+z^3t)+1) \epm.
$$
Let $\{e_1, e_2, e_3\}$ be the standard basis on $\R^{3}$, and $V_1=\R v_1$,  $V_2=\R e_2 \oplus \R e_3$. Then $V_1\subset V_1^\perp$ and $V_2^\perp\subset V_2$. Let $\a_1 \neq \a_2 \in \R$. Set $\ti p_1 =E^{-1}(x, t, \a_1^2)e_1$, $\ti p_2=E^{-1}(x, t, \a_2^2)e_2$, $\ti p_3=E^{-1}(x, t, \a_2^2)e_3$. And 
$$\ti V_1=\R\ti p_1, \quad \ti V_2=\R\ti p_2\oplus\R\ti p_3.$$ 
The projection $\ti \pi$ of $\R^{2, 1}$ onto $\ti V_1$ with respect to $\R^{2, 1}=\ti V_1 \oplus \ti V_2$ is 
$$\ti \pi=\bpm \ti p_1 & 0 & 0\epm \bpm \ti p_1 & \ti p_2 & \ti p_3\epm^{-1}.$$
To simplify the result, we introduce some notation:
\beq\label{cv}
c(x, t, \a)=\cosh(\a x+\a^3 t), \quad s(x, t, \a)=\sinh (\a x+\a^3 t).
\eeq
Let $D$ denote the determinant of $\bpm \ti p_1 & \ti p_2 & \ti p_3\epm$. A direct computation implies that 
$$D=\frac{1}{4}(c(\a_1)+1)(c(\a_2)+1)-\frac{\a_2}{2 \a_1}s(\a_1)s(\a_2)+\frac{\a_2^2}{4 \a_1^2}(c(\a_1)-1)(c(\a_2)-1).$$
Apply  BT for \eqref{ah} to get a new solution $\ti q(x,t)=\bpm  q_1 & q_2 & 0 \\ 0 & 0 & q_2 \\ 0 & 0 & -q_1\epm$ of \eqref{ah} (with $j=2$), where
\begin{align*}
& q_1=\frac{(\a_1^2-\a_2^2)(\a_2(c(\a_1)-1)s(\a_2)-\a_1s(\a_1)(c(\a_2)+1))}{4\a_1^2 D}, \\
& q_2=\frac{(\a_1^2-\a_2^2)(\a_2^2(c(\a_1)-1)(c(\a_2)-1)-\a_1^2(c(\a_1)+1)(c(\a_2)+1))}{8\a_1^2D}.
\end{align*}
Here we use $c(\a)$, and $s(\a)$ to denote functions $c(\cdot, \cdot, \a)$ and $s(\cdot, \cdot, \a)$ respectively. 
So the new solution is $\ti u= y(e_{12}+ e_{21})$ for \eqref{ai}, where  
$y=q_2+\frac{1}{2}q_1^2+(q_1)_x$. And the corresponding curve flow solution for \eqref{bx} is:
$$\ti \g=\frac{\a_1^2-\a_2^2}{2 \a_1^2 D} \bpm 1 & 0 & 0 \\ x & 1 & 0 \\ \frac{x^2}{2} & x & 1 \epm
\bpm \frac{2 \a_1^2 D}{\a_1^2-\a_2^2} - c(\a_1)-c(\a_2)\\
\frac{1}{\a_1}s(\a_1)(c(\a_2)+1)-\frac{1}{\a_2}(c(\a_1)-1)s(\a_2)\\
\frac{\a_1^2-\a_2^2}{\a_1^2\a_2^2}(c(\a_1)-1)(c(\a_2)+1)
\epm$$
In particular, let $\a_2=0$. We get the following smooth solution of \eqref{bx},
$$\ti \g(x,t)=\bpm 1-\frac{\a s(x,t,\a)}{c(x,t,\a)+1}x-\frac{\a^2(c(x,t,\a)-1)}{4(c(x,t,\a)+1)}x^2  \\
\frac{2s(x,t,\a)}{\a(c(x,t,\a)+1)}-\frac{c(x,t,\a)-1}{c(x,t\a)+1}x \\
\frac{2(c(x,t,\a)-1)}{\a^2(c(x,t,\a)+1)}
\epm, $$
where $\a\in \R$, and  $c(x,t,\a)$ $s(x,t,\a)$ are defined by \eqref{cv}.
The isotropic curvature of $\ti \g$ is the one-soliton of the KdV,
$$\ti u=-\a_1^2\sech^2(\frac{a_1}{2}x+\frac{\a_1^3}{2}t).$$
\eeg

Next we give a Permutability formula for BTs of \eqref{ah}. First we need some Lemmas.

\blem\label{bl} 
Let $\a_i, \b_i, i=1, 2$ be four distinct constants in $\R \bh \{0\}$, and $V_{i}, W_{i}$ linear subspaces of $\R^{n+1, n}$ for $i=1, 2$ such that 
\begin{align*}
& \R^{n+1, n}=V_1 \oplus V_2, \quad V_1 \subset V_1^{\perp}, \quad V_2^{\perp} \subset V_2, \\
& \R^{n+1, n}=W_1 \oplus W_2, \quad W_1 \subset W_1^{\perp}, \quad W_2^{\perp} \subset W_2.
\end{align*}
Let $\pi_1$ be the projection of $\R^{n+1, n}$ onto $V_1$ along $V_2$, and $\pi_2$ the projection of $\R^{n+1, n}$ onto $W_1$ along $W_2$. Set
$$\bca
\hat W_1=g_{\b_1, \b_2, \pi_2}(\a_1)(V_1),  \\
\hat W_2=g_{\b_1, \b_2, \pi_2}(\a_2)(V_2).
\eca
$$
Then
\ben 
\item[(i)] $\hat W_1 \subset \hat W_1^{\perp},  \hat W_2^{\perp} \subset \hat W_2$,
\item[(ii)] $\hat W_2^{\perp}=g_{\b_1, \b_2, \pi_2}(\a_2)(V_2^{\perp})$.
\een
\elem

\begin{proof} Note that $g_{\b_1, \b_2, \pi_2}(\a_i) \in O(n+1, n)$, for $i=1, 2$. So for any $X, Y \in \R^{n+1, n}$, 
$$\li g_{\b_1, \b_2, \pi_2}(\a_i) , g_{\b_1, \b_2, \pi_2}(\a_i) Y\ri=\li X, Y \ri, \quad i=1, 2.$$
This proves (i). To prove (ii), note that $X \in \hat W_2^{\perp}$ if and only if 
$$\li X, g_{\b_1, \b_2, \pi_2}(\a_2)Y\ri=\li g_{\b_1, \b_2, \pi_2}(\a_i)^{-1}X, Y \ri=0, \quad \forall \ Y \in V_2.$$
Hence $g_{\b_1, \b_2, \pi_2}(\a_i)^{-1}X \in V_2^{\perp}$. This proves (ii).
\end{proof}

\bprop\label{bm} Let $\a_1, \a_2, \b_1, \b_2$ be four distinct non-zero real numbers, and $\pi_i, \tau_i$ projections of $\R^{n+1, n}$  satisfying $\pi_i\pi_i^\sharp= \pi_i^\sharp\pi_i = \tau_i \tau_i^\sharp =\tau_i^\sharp \tau_i = {\bf 0}$ for $1\leq i\leq 2$. Set 
\begin{align*}
&V_1= \Im(\pi_1), V_2= \Ker(\pi_1), W_1= \Im(\pi_2), W_2=\Ker(\pi_2), \\
& \hat V_1= \Im(\hat \tau_2), \hat V_2=\Ker(\hat\tau_2),  \hat W_1= \Im(\tau_1), \hat W_2= \Ker(\hat\tau_1).
\end{align*}
 Then 
\beq\label{gf}
g_{\b_1, \b_2, \tau_2}g_{\a_1, \a_2, \pi_1}=g_{\a_1, \a_2, \tau_1}g_{\b_1, \b_2, \pi_2}
\eeq
if and only if 
\beq\label{gm}
\bca
\hat V_1=g_{\a_1, \a_2, \tau_1}(\b_2)W_1,  \quad
\hat V_2=g_{\a_1, \a_2, \tau_1}(\b_1)W_2,\\
\hat W_1=g_{\b_1, \b_2, \pi_2}(\a_1)(V_1),  \quad
\hat W_2=g_{\b_1, \b_2, \pi_2}(\a_2)(V_2).
\eca
\eeq
\eprop

\begin{proof} We first prove \eqref{gm} is sufficient.  The proof of Lemma \ref{bl} implies that $\hat V_1 \subset \hat V_1^{\perp},  \hat V_2^{\perp} \subset \hat V_2$. Since $g_{\a_1, \a_2, \pi_1}^{-1}=g_{\a_1, \a_2, \pi_1^{\sharp}}$, it suffices to prove that 
\beq\label{bn}
g_{\b_1, \b_2, \pi_2}g_{\a_1, \a_2, \pi_1^{\sharp}}=g_{\a_1, \a_2, \tau_1^{\sharp}}g_{\b_1,\b_2, \tau_2}.
\eeq

It is equivalent to prove that the residues at $\a_1, \b_2$ on both sides of \eqref{bn} are equal, i.e.,
$$
\bca
g_{\b_1, \b_2, \pi_2}(\a_1)\pi_1=\tau_1g_{\b_1, \b_2, \tau_2^{\sharp}}(\a_1), \\ g_{\b_1, \b_2, \pi_2}(\a_2)\pi_1^{\sharp}=\tau_1^{\sharp}g_{\b_1, \b_2, \tau_2^{\sharp}}(\a_2),
\eca
\bca
\pi_2^{\sharp}g_{\a_1, \a_2, \pi_1}(\b_1)=g_{\a_1, \a_2, \tau_1^{\sharp}}(\b_1)\tau_2^{\sharp}, \\
\pi_2g_{\a_1, \a_2, \pi_1}(\b_2)=g_{\a_1, \a_2, \tau_1^{\sharp}}(\b_2)\tau_2.
\eca
$$
These are true because we have
\begin{align*}
& \Im(\tau_1)=g_{\b_1, \b_2, \pi_2}(\a_1)\Im(\pi_1), \quad
\Im(I-\tau_1)=g_{\b_1, \b_2, \pi_2}(\a_2)\Im(I-\pi_1), \\
& \Im(\tau_2)=g_{\a_1, \a_2, \tau_1}(\b_2)\Im(\pi_2), \quad
\Im(I-\tau_2)=g_{\a_1, \a_2, \tau_1}(\b_1)\Im(I-\pi_2).
\end{align*}

The computation given for the sufficient part also proves necessary part. 
\end{proof}

\bthm\label{ev} {\bf [Permutability]} \ 

\ni Let $q$ be a solution of  \eqref{ah}, and $\a_i, \b_i, \pi_i, \tau_i$  as  in Proposition \ref{bm} for  $i=1, 2$ satisfying \eqref{gf}.  Let
\begin{align*}
q_1&=g_{\a_1, \a_2, \pi_1} \sharp q = q+ (\a_1-\a_2)[\b, \ti \pi_1-\ti\pi_1^\sharp], \\
q_2 &=g_{\b_1, \b_2, \pi_2} \sharp q = q+ (\b_1-\b_2)[\b, \ti \pi_2-\ti\pi_2^\sharp],
\end{align*}
 $q_{12}=g_{\b_1, \b_2, \tau_2}\sharp(g_{\a_1, \a_2, \pi_1}\sharp q)$ and $q_{21}= g_{\a_1, \a_2,\tau_1}\sharp(g_{\b_1,\b_2,\pi_2}\sharp q)$ solutions of \eqref{ah} constructed from B\"acklund transformations (Theorem \ref{bea}). 
Let 
$\ti \tau_1$ and $\ti \tau_2$ be the projections of $\R^{n+1, n}$ with
\begin{align}
& \Im(\ti\tau_1)=g_{\b_1, \b_2, \ti\pi_2}(\a_1)\Im(\ti \pi_1), \quad
\Ker(\ti\tau_1)=g_{\b_1,\ti \b_2, \ti\pi_2}(\a_2)\Ker(\ti \pi_1), \label{gp}\\
&
\Im(\ti\tau_2)=g_{\a_1, \a_2, \ti \pi_1}(\b_2)\Im(\ti \pi_2), \quad
\Ker(\ti\tau_2)=g_{\a_1, \a_2, \ti \pi_1}(\b_1)\Ker(\ti \pi_2).  \label{gq}
\end{align} 
Then 
\beq\label{gr}
q_{12}= q_{21} = q_1+ (\b_1-\b_2)[\b, \ti\tau_2-\ti \tau_2^\sharp] = q_2+ (\a_1-\a_2)[\b, \ti\tau_1-\ti\tau_1^\sharp],
\eeq
where $\b$ is as in \eqref{db}. 
\ethm

\begin{proof} Let $F$ be a frame of the solution $q$ of \eqref{ah}. Theorem \ref{bea} implies that 
\begin{align*}
F_1& = g_{\a_1, \a_2, \pi_1}F g_{\a_1, \a_2, \ti\pi_1}^{-1},\\
F_2& = g_{\b_1, \b_2, \pi_2}F g_{\b_1, \b_2, \ti\pi_2}^{-1},
\end{align*}
are frames of $q_1$ and $q_2$ respectively.  Apply Theorem \ref{bea} to $q_2$ and $q_1$ to see that there are projections $\ti\tau_1(x,t)$ and $\ti\tau_2(x,t)$ such that
\begin{align*}
F_{12} &=g_{\b_1, \b_2, \tau_2}g_{\a_1, \a_2, \pi_1} F g^{-1}_{\a_1, \a_2, \ti \pi_1}g^{-1}_{\b_1, \b_2, \ti\tau_2},\\
F_{21} & = g_{\a_1, \a_2, \tau_1} g_{\b_1, \b_2, \pi_2} Fg^{-1}_{\b_1, \b_2, \ti\pi_2} g^{-1}_{\a_1, \a_2, \ti\tau_1}
\end{align*}
are frames of $q_{12}$ and $q_{21}$ respectively.  Let $f=g_{\b_1, \b_2, \tau_2}g_{\a_1, \a_2, \pi_1}$. By assumption, 
$f=  g_{\a_1, \a_2, \tau_1} g_{\b_1, \b_2, \pi_2}$. So we obtain 
$$fF= F_{12} g_{\b_1, \b_2, \ti\tau_2}g_{\a_1, \a_2, \ti\pi_1} = F_{21}g_{\a_1, \a_2, \ti\tau_1} g_{\b_1,\b_2, \ti \pi_2}.$$
This gives two factorizations of $fF$ as the product of elements in $(\bbn1)_+$ and $(\bbn1)_-$. Since the factorization of $fF$ in $(\bbn1)_+(\bbn1)_-$ is unique, we get $F_{12}= F_{21}$ and
\beq\label{go} g_{\b_1,\b_2, \ti\tau_2} g_{\a_1, \a_2, \ti\pi_1} = g_{\a_1, \a_2, \ti\tau_1} g_{\b_1,\b_2, \ti \pi_2}. 
\eeq
It follows from \eqref{go} and Proposition \ref{bm} that $\ti\tau_i$ satisfies \eqref{gp} and \eqref{gq}.  Since $F_{12}=F_{21}$, we have $q_{12}= q_{21}$. Formula \eqref{gr} follows from Theorem \ref{bea}. 
\end{proof}

\brem \label{ey}\hfil

\ben
\item[(i)] The solution $q_{12}$ given Theorem \ref{ev} is an algebraic function of $q, \ti \pi_1, \ti\pi_2$. 
\item[(ii)] We apply Theorem \ref{bea} to the trivial solution $q=0$ to get $k$ $1$-soliton solutions $q_i$ of \eqref{ah} and their frames $F_i$ for $1\leq i\leq k$. Apply the Permutability Theorem \ref{ev} to construct 2-soliton solutions $q_{ij}$ of \eqref{ah}. Apply Theorem \ref{ev} to $q_2$, $q_{12}$ and $q_{23}$ to get the 3-soliton solution $q_{123}$. Continue this way to get explicit formulas for $k$-soliton solutions of \eqref{ah} and their frames algebraically from one-solitons. Suppose $\ti F$ is a frame of a $k$-soliton solution  $\ti q$ of \eqref{ah}, and $\ti \D= D(\ti q)$, where $D$ is the operator given in Definition \ref{aba}. Then $\ti E= \ti F\ti \D^{-1}$ is a frame of the solution $\ti u=\ti \D \ast \ti q$ of the $(2j-1)$-th $\bn1$-KdV flow \eqref{ai} and $\ti \g(x,t)= \ti E(x,t,0)e_1$ is an explicit solution of the isotropic curve flow \eqref{cb}. 
\een
\erem

\bs

\section{B\"acklund transformations for the $\an2$-KdV flows}\label{hg}

We proceed the same way as for the $\bn 1$-KdV flows in section \ref{hd} to construct BTs for isotropic curve flows of A-type on $\cm_{n+1, n}$.

Let $\ban2$ denote the group of smooth maps $f:S^1\to SL(2n+1,\C)$ satisfying 
\beq\label{bf}
\overline{f(\bar \l)}= f(\l),  \quad f(-\l)^t C_n f(\l)= C_n,
\eeq 
and $(\ban2)_+$  the subgroup of $f\in \ban2$ that can be extended to a holomorphic map on $\C$, and $(\ban2)_-$ the subgroup of $f\in \ban2$ that is the boundary value of a holomorphic map $\ti f$ defined on $\e^{-1}\leq |\l |\leq \infty$ such that $\ti f(\infty)=\I$ for some small $\e>0$. Then the Lie algebras of $\ban2$ and $(\ban2)_\pm$ are $\an2$ and $(\an2)_\pm$ respectively.  

\bdefn If $V$ is a linear subspace of $\R^{n+1, n}$ such that $\R^{n+1, n}= V\oplus V^\perp$, then we call the projection of $\R^{n+1, n}$ onto $V$ along $V^\perp$ an {\it $O(n+1, n)$-projection\/}. 
\edefn

The proofs of the following two Propositions are straight forward.

\bprop \label{cw}\

\ben
\item[(a)] $\pi$ is an $O(n+1, n)$-projection if and only if  $\pi^2=\pi=\pi^\sharp$.
\item[(b)] If $V$ is a linear subspace of $\R^{n+1, n}$, then the restriction of $\li\, ,\ri$ to $V$ is non-degenerate if and only if  $\R^{n+1, n}= V\oplus V^\perp$.
\een
\eprop

\bprop\label{cx}
Let $\pi$ be an $O(n+1, n)$-projection of $\R^{n+1, n}$ onto $V$ along $V^\perp$, and $\alpha \in \R \backslash \{0\}$. Set
\beq\label{bfa}
g_{\a, \pi}(\l)=I+\frac{2\a}{\l - \a}(I-\pi).
\eeq
Then $g_{\a, \pi} \in (\an2)_-$.
\eprop

\bthm \label{bg}
Let $F(x,t,\cdot) \in (\hat A_{2n}^{(2)})_+$ be a frame of a
solution $q$ of \eqref{ak}, $\pi$ be an $O(n+1, n)$-projection onto $V$ along $V^\perp$, $\a\in \R\bh \{0\}$ a constant, and $g_{\a, \pi}$ defined as in \eqref{bfa}. Set
$\ti{V}(x, t)=F(x,t,\a)^{-1}(V)$, $\ti\pi(x,t)$ the $O(n+1, n)$-projection onto $\ti V(x,t)$ and $\ti F(x, t, \l)=g_{\a, \pi}(\l)F(x, t, \l)g_{\a, \ti \pi(x, t)}(\l)^{-1}$. Then 
\ben
\item[(1)] $\R^{n+1, n}= \ti V(x,t)\oplus \ti V(x,t)^\perp$,
\item[(2)] $\ti q=q+2\a[\ti \pi, e_{1, 2n+1}]$ is a solution of \eqref{ak} and $\ti F$ is a frame of $\ti q$ (We will use $g_{\a,\pi}\bu q$ to denote the new solution $\ti q$).
\een
\ethm

\begin{proof} Since $F(x,t,\cdot)$ satisfies \eqref{bf} and $\a\in \R$, $F(x,t,\a)\in O(n+1, n)$. Hence the restriction to  $\ti V(x,t)$ is non-degenerate..  This proves statement (1).  

By Proposition \ref{cx}, $g_{\a, \ti\pi(x,t)}\in (\an2)_-$. So to prove (2), it suffices to prove that $\ti F(x,t,\l)$ is holomorphic in $\l$.  To prove this, we only need to show that the residues of $\ti F(x, t, \l)$ at $\pm \a$ are zero. Note that the residue of $\ti F(x, t, \l)$ at $\a$ is 
$$2\a (I- \pi)F(x, t, \a)\ti{\pi}.$$
Since $\Im (\ti \pi(x,t))=\ti V(x,t)$ and $F(x, t, \a)\ti V(x,t)=V$, we have  
$$F(x, t, \a)\Im(\ti V(x,t))=V.$$ 
Therefore the residue of $\ti F(x, t, \l)$ at $\l=\a$ is zero. 

The residue of $\ti F(x, t, \l)$ at $-\a$ is 
$$ 2\a\pi F(x,t, -\a)(I -\ti \pi).$$
For any $\eta \in \ti V(x,t)^\perp$, $\li \eta, \hat v\ri=0$, since $F(x, t, \cdot) \in (\an 2)_+$, we have 
$$\li F(x, t, -\a)\eta, V\ri = \li F(x, t, -\a)\eta, F(x, t, \a)\ti V(x,t)\ri=\li \eta, \ti V(x,t)\ri=0.$$ 
Hence $\pi F(x,t, -a)(I -\ti \pi)\eta=0$. This proves the residue of $\ti F(x, t, \l)$ at $\l=-\a$ is zero.  This proves $\ti F(x,t,\l)$ is holomorphic for $\l\in \C$
\end{proof}

\bthm\label{dp} Let $E(x,t,\cdot)\in (\ban2)_+$ be a frame of a solution $u$ of \eqref{cf}, i.e., $u_t= [\p_x+b+u, S_{2j-1,0}(u)-\ti\eta_j(u)]$. Let $\D(x,t)\in N_n^+$ satisfying $\D_t\D^{-1}= \ti\eta_j(u)$.  Then
\ben
\item[(i)]  $q=\D^{-1}\ast u$. is a solution of \eqref{ak} and $F= E\D$ a frame for $q$.
\item[(ii)]  Let $\a, \pi, \ti \pi$ be as in Theorem \ref{bg}, and $\ti q= g_{\a, \ti\pi}\bu q$. Then 
$$\ti E= g_{\a, \pi}E\D g_{\a, \ti\pi}^{-1}\ti \D^{-1}$$
is a frame of a new solution $\ti u= \ti \D\ast(g_{\a, \pi}\bu (\D \ast u))$ of \eqref{cf}.
\een
 \ethm

\bthm\label{ew} Let $\g: \R^2 \rightarrow \cm_{n+1, n}$ be a solution of the $(2j-1)$-th isotropic curve flow  \eqref{ca} of  A-type, $g(\cdot, t)$ the isotropic moving frame along $\g(\cdot, t)$, and $u(\cdot, t)=(g^{-1}g_x-b)$ the solution of the $(2j-1)$-th $\an2$-KdV flow \eqref{cf} as in Theorem \ref{an}. Let $\a, \pi, \ti \pi$, $\D$, $\ti E$ and $\ti u$ be as in Theorem \ref{bg}. Then $\ti\g(x,t)=\ti E(x,t, 0)e_1$ is is a solution of \eqref{ca} and $\ti u$ is its  isotropic curvature, where  $e_1=(1, 0, \cdots, 0)^t$.
\ethm

\bcor Let $\g, g, \a_1, \a_2, \pi, \ti \pi, \D$ be as in Theorem \ref{ew}. If $\a_1\a_2\not=0$, then 
$$ \ti{\g}(x, t)=g(x, t)\D(x,t) (2\ti\pi(x,t)-I_{2n+1})e_1$$
is a solution of \eqref{ca}.
\ecor

\beg {\bf $1$-soliton of the isotropic curve flow of A-type on $\cm_{2,1}$ of A-type} \ 

Note that $\g=(1, x, \frac{x}{2})^t \in \cm_{2,1}$ is a solution of \eqref{al} with isotropic curvature $q=0$. Set $\l =z^3$, 
\begin{align*}
& D(z)= \diag(1, z, z^2), \quad  \Xi= (\a^{(i-1)(j-1)})_{3 \times 3}, \quad   \sigma=e_{12}+e_{23}+e_{31}, \\
& A_i(x, t, z)=\exp( \a^{i-1} zx+ (\a^{i-1} z)^5 t), \quad i=1, 2, 3,\\ 
&(m_1(x,t, z), m_2(x, t, z), m_3(x,t, z))
= (e^{A_1}, \ldots, e^{A_3})\Xi.
\end{align*}
Then the extended frame $E(x, t, z^3)$ of $q=0$ is 
\begin{align*}
E(x, t, z^3)=\frac{1}{3}D(z)^{-1}\bpm m_1(x, t, z) & m_2(x, t, z) & m_3(x, t, z) \\ m_3(x, t, z) & m_1(x, t, z) & m_2(x, t, z) \\ m_2(x, t, z) & m_3(x, t, z) & m_1(x, t, z) \epm D(z). 
\end{align*}
Apply Theorem \ref{bg} by choosing $V=\R v$ for some constant $v \in \R^{2, 1}$ with $\li v, v\ri \neq 0$. Then we have
$$\ti v=E(-x, -t, k^3)v, \quad \ti v_i= \frac{1}{3}(e^{A_1}, e^{A_2}, e^{A_3})\Xi\sigma^{i-1} D(k)(v_1, v_2, v_3)^t.$$

In particular, we can choose $v$ such that 
\begin{align*}
\ti v= (e^{A_1}+\a^2 e^{A_2}, \frac{1}{k}(e^{A_1}+\a e^{A_2}), \frac{1}{k^2}(e^{A_1}+e^{A_2}) )^t.
\end{align*}
Apply BT for \eqref{ak} (Theorem \ref{bg}), we get a new solution $\ti q=\ti q_1 (e_{11}-e_{33})+ \ti q_2 (e_{12}+e_{23})$, where
\begin{align*}
& \ti q_1= 2k\frac{(e^{A_1}+e^{A_2})^2}{4\a e^{A_1+A_2}-e^{2A_1}-\a^2e^{2A_2}}\\
& \ti q_2=-2k^2\frac{e^{2A_1}+\a e^{2A_2}+(1+\a)e^{A_1+A_2}}{4\a e^{A_1+A_2}-e^{2A_1}-\a^2e^{2A_2}}.
\end{align*}
Hence the new solution $\ti u= y(e_{12}+ e_{21})$ of \eqref{al}, where  
$$y=\ti q_2+\frac{1}{2}\ti q_1^2+(\ti q_1)_x.$$
Let $M=4\a e^{A_1+A_2}-e^{2A_1}-\a^2e^{2A_2}$, the corresponding curve flow solution for \eqref{al} is 
$$\ti \g=\bpm 1 & 0 & 0 \\ x & 1 & 0 \\ \frac{x^2}{2} & x & 1 \epm \bpm \frac{2}{M}(e^{2A_1}+(1+\a^2)e^{A_1+A_2}+\a^2e^{2A_2}+\frac{1}{2k^3}\ti q^2 \\
\frac{1}{2k^3}\ti q_2 \\ -\frac{1}{2k^3}\ti q_1\epm$$
\eeg

Next we derive the permutability formula for BTs of \eqref{ak}. First we need a Proposition.

\bprop \label{bi} Let $\a_1, \a_2 \in \R\backslash \{0\}$, and $\pi_i$ an $O(n+1, n)$-projection of $\R^{n+1, n}$ for $i=1, 2$.  If $|\a_1| \neq |\a_2|$, then 
$$\phi=\a_1-\a_2+2\a_2\pi_2-2\a_1 \pi_1$$
is invertible. Moreover,
\ben
\item[(i)] $\tau_i=\phi \pi_i \phi^{-1}$ is an $O(n+1, n)$-projection for $i=1, 2$, 
\item[(ii)]
$g_{\a_2, \tau_2}g_{\a_1, \pi_1}=g_{\a_1, \tau_1}g_{\a_2, \pi_2}$ if and only if $\tau_i=\phi \pi_i \phi^{-1}$ for $i=1,2$.
\een
\eprop

\begin{proof} Since $\pi_i^2=\pi_i$, the eigenvalues of $\pi$ are $0, \pm 1$. Consider the eigenvalue of $\phi$, if $\pi_1=\pi_2$, then the eigenvalue of $\phi$ is $\a_2-\a_1$. If $\pi_1 \neq \pi_2$, then the possible eigenvalues for $\pi$ are $\pm(\a_1+\a_2)$. Therefore, as long as $|\a_1| \neq |a_2|$, $\phi$ is invertible.

From the definition of $\tau_i$, we have $\tau_i^2=\tau_i$ for $i=1, 2$. To prove (i), we only need to show that $\tau_i^t C_n=C_n \tau_i$. It is equivalent to prove that 
\beq\label{bh}
\pi_i^t\phi^tC_n \phi=\phi^t C_n \phi \pi_i, \quad i=1, 2.
\eeq
It follows from direct computation and the fact that both $\pi_1$ and $\pi_2$ are $O(n+1, n)$-projections, we have
$$\phi^tC_n \phi=(\a_1-\a_2)^2C_n+4\a_1\a_2C_n(\pi_1+\pi_2-\pi_2\pi_1-\pi_1\pi_2).$$
Hence 
\begin{align*}
\pi_1^t\phi^tC_n\phi & =(\a_1-\a_2)^2C_n\pi_1+4\a_1\a_2C_n(\pi_1+\pi_1\pi_2-\pi_1\pi_2\pi_1-\pi_1\pi_2) \\
\quad & = \phi^tC_n\phi \pi_1
\end{align*}
Similarly, $\pi_2^t\phi^tC_n\phi=\phi^tC_n\phi \pi_2$. 

To prove (ii), let
\begin{align*}
& Y_1=\a_1-2\a_1\pi_1, \quad Y_2=\a_2-2\a_2\tau_2, \\
& Z_1=\a_1-2\a_1\tau_1, \quad Z_2=\a_2-2\a_2\pi_2.
\end{align*}
We claim that $(\l+Y_2)(\l+Y_1)=(\l+Z_1)(\l+Z_2)$. Compare coefficients as an expansion of $\l$ to get
$$
\bca
Y_1+Y_2=Z_1+Z_2, \\
Y_2Y_1=Z_1Z_2.
\eca
$$
So we get
$$
Z_1=(Y_1-Z_2)Y_1(Y_1-Z_2)^{-1}, \quad Y_2=(Y_1-Z_2)Z_2(Y_1-Z_2)^{-1}.
 $$
 Hence $\tau_i=\phi \pi_i \phi^{-1}, i =1, 2$.
\end{proof}

\bthm \label{ez} {\bf [Permutability]} \

\ni Let $\a_1, \a_2 \in \R\backslash\{0\}$ with $|\a_1| \neq |\a_2|$, and $\pi_i$ an $O(n+1, n)$-projection for $i=1, 2$. Let $q$ be a solution of \eqref{ak},  
$q_i=g_{\a_i, \pi_i}\bu q= q+ 2\a_i[\ti\pi_i, e_{1, 2n+1}]$, and 
$$q_{12}= g_{\a_2, \tau_2}\bu(g_{\a_1, \pi_1}\bu q), \quad q_{21}=g_{\a_1, \tau_1}\bu(g_{\a_2, \pi_2}\bu q)$$
 solutions of \eqref{ak} obtained from B\"acklund transformations. Set
\begin{align*}
& \phi=\a_1-\a_2+2\a_2\pi_2-2\a_i\pi_1, 
\quad \ti\phi=\a_1-\a_2+2\a_2\ti\pi_2-2\a_i\ti\pi_1, \\
& \tau_i=\phi \pi_i \phi^{-1}, \quad \ti\tau_i=\phi \ti\pi_i \phi^{-1}, \quad i=1, 2.
\end{align*}
Then 
$$q_{12}= q_{21} = q_1+2\a_2[\ti \tau_2, e_{1, 2n+1}] = q_2+ 2\a_1[\ti \tau_1, e_{1, 2n+1}].$$
\ethm

Similarly, we can use Theorems \ref{bg} and \ref{ew} to construct explicit $k$-soliton solutions for isotropic curve flows of A-type. 



\begin{thebibliography}{99}

\bibitem{CIM13}Calini, A., Ivey, T., Mar{\'{\i}} Beffa,G., \emph{Integrable flows for starlike curves in centroaffine space}, SIGMA Symmetry Integrability Geom. Methods Appl. \textbf{9} (2013), Paper 022, 21 pp.

\bibitem{YC92}
Cheng, Y., \emph{Constraints of the Kadmotesev-Petviashvili hierarchy}, J. Math. Phys. \textbf{33(11)}(1992), 3774--3782.

\bibitem{Dic03}Dickey, L. A., \emph{Soliton equations and Hamiltonian systems}, second edition, Advanced Series in Mathematical Physics \textbf{26} (2003), World Scientific Publishing Co. Inc., River Edge, NJ.

\bibitem{DS84}Drinfel'd, V.G., Sokolov, V.V., \emph{Lie algebras and equations of Korteweg-de Vries type},  (Russian) Current problems in mathematics, \textbf{24} (1984),  81--180, Itogi Nauki i Tekhniki, Akad. Nauk SSSR, Vsesoyuz. Inst. Nauchn. i Tekhn. Inform., Moscow.

\bibitem{K80}Kaup, D. J., \emph{On the inverse scattering problem for cubic eigenvalue problems of the class {$\psi _{xxx}+6Q\psi _{x}+6R\psi =\lambda \psi $}}, Stud. Appl. Math. \textbf{62(3)} (1980),189--216.

\bibitem{Ku84}Kupershmidt, B.A., \emph{A super Korteweg-de Vries equation: an integrable system}, Phys. Lett. A \textbf{102(5--6)}  (1984), 213--215.

\bibitem{UP95}Pinkall, U., \emph{Hamiltonian flows on the space of star-shaped curves},  Results Math. \textbf{27(3-4)} (1995), 328--332.

\bibitem{TU00} Terng, C.L., Uhlenbeck, K., \emph{B\"acklund transformations and loop group actions}, Comm. Pure Appl. Math. \textbf{53(1)} (2000), 1--75.

\bibitem{TU11}Terng, C.L., Uhlenbeck, K., \emph{The $n \times n$ KdV hierarchy}, JFPTA \textbf{10} (2011), 37--61.

\bibitem{TWa}Terng, C.L., Wu, Z., \emph{Central affine curve flow on the plane}, J. Fixed Point Theory Appl., Mme Choquet-Bruhat Fastschrift, \textbf{14} (2013), 375--396. 

\bibitem{TWb}
Terng, C.L., Wu, Z., \emph{N-dimension central affine curve flows}, arXiv:1411.2725v1. 
\end{thebibliography}
\end{document}